\documentstyle{amsppt}
\magnification=\magstephalf \pagewidth{15 truecm} \pageheight{23
truecm} \NoRunningHeads \nologo \loadbold \refstyle{A}
\NoBlackBoxes \TagsOnRight \font\medium=cmbx12 scaled \magstephalf

\def\d{\delta}
\def\f{\frac}

\def\l{\lambda}
\def\L{\Lambda}
\def\O{\Omega}
\def\ri{\right}

\def\t{\tilde}
\def\th{\theta}
\def\ve{\varepsilon}

\def\re{\par\hang\textindent}

\def\spb{\smallpagebreak}
\def\mpb{\vskip 0.5truecm}

\topmatter
\title
{\medium {The symplectic geometry of a new kind of Siegel upper
half space of order 2 (I)}}
\endtitle
\author
Tianqin Wang, Tianze Wang, and Hongwen Lu
\endauthor

\address
\endaddress

\abstract {In this paper, we introduce a new kind of Siegel upper
half space and consider the symplectic geometry on it explicitly
under the action of the group of all holomorphic transformations
of it. The results and methods will form a basis for our number
theoretic applications later.}
\endabstract

\thanks
Project Supported by the National Natural Science Foundation of
China, Grant Number 11471112.
\endthanks
\endtopmatter

\document

{\bf Keywords:}\,\,\, {Siegel upper half space; symplectic
geometry; group action.}

\mpb {\bf Mathematics Subject Classification 2010:}
11F55,\,\,\,11F46,\,\,\,11F50.

\vskip 0.6cm

\head 1. Introduction and notations
\endhead

Let $\Bbb{R}$ denote the field of real numbers, and $\Bbb{C}$ the
field of complex numbers by convention. For any positive integers
$m$ and $n$, let $\Bbb{C}^{(m,n)}$ denote the set of all $m\times
n$ matrices $Z$ with entries in $\Bbb{C}$. For any
$Z\in\Bbb{C}^{(m,n)}$ we use $^{t}Z$ to denote its transpose,
which is a $n\times m$ matrix in $\Bbb{C}^{(n,m)}$. For any
positive integer $n$, the well known so-called Siegel upper half
space $\Bbb{H}_{n}$ of order $n$ is, by definition,
$$
\Bbb{H}_n:=\left\{Z\in\Bbb{C}^{(n,n)}\mid\  ^{t}Z=Z,\ {\text Im}
Z>0\right\},
$$
where ${\text Im} Z$ (resp. ${\text Re} Z$) denotes the imaginary
(resp. real) part of $Z$, i.e., if $Z=(z_{ij})$ then ${\text Im}
Z=({\text Im} z_{ij}),\ {\text Re} Z=({\text Re} z_{ij})$, and
${\text Im} Z>0$ means that ${\text Im} Z$ is positive definite.
Here, and throughout this paper, the symbol $":="$ is used to
indicate that the right hand side of an equality is the definition
of the left. Clearly, if $n=1$, then $\Bbb{H}_{n}$ is reduced to
the classical Poinc{\'a}re[P] upper half plane
$\Bbb{H}=\Bbb{H}_1$.

The symplectic geometry on $\Bbb{H}$ by the action of the group of
linear fractional transformations, or, by the order $2$ symplectic
group $Sp(2,\Bbb{R})=SL_2(\Bbb{R})$, possesses an almost full
satisfactory understanding nowadays. And this understanding has
led to many great applications both in mathematics itself and the
spread area of other branches of subject. One of the most
fascinating examples is the basic role it plays to the theory of
classical automorphic forms and the applications in the field of
number theory.

For $n>1$, the symplectic geometry on $\Bbb{H}_n$, acted by the
order $2n$ symplectic group $Sp(2n,\Bbb{R})$, were studied
systematically by Siegel[S2] in 1936 for the first time, and lots
of very important essential results were established. Besides many
other important applications it leads to, it then becomes the
concrete basis for the theory of Siegel modular forms and Jacobi
forms, see, e.g. [AZ], [DI], [EZ], [M], [S3], [Sk3], [Sk4], [SZ1]
and [SZ2]. And these subjects nowadays are becoming more and more
important, active and fertile fields of mathematics.

In this series of papers, we will introduce a new kind of Siegel
upper half space $\hat\Bbb H_2$ of order 2, see (3) below; and
give a relatively systematic argument for the symplectic geometry
on it, in the light of Siegel's classical work[S2]. This first
paper will focus only on the relatively pure geometric part of our
work. Although one might notice that our new object is not
irreducible by Cartan[C] and the geometry might be viewed as a
topological product of two Poinc{\'a}re upper half planes $\Bbb H$
in essence as we proceed, it should be pointed out in advance that
there do exist many important interesting and useful results
valuable for researching at least from the explicit point of view.
Furthermore, besides the independent importance just from the
geometric point of view, one will find that the geometry in this
paper will naturally become the essential basis for our number
theoretic applications that follow in this series later. In other
words, the principal novelties of this paper can be realized
mainly from the following two aspects of view: one is a systematic
and explicit formulation of the principal geometry of the new kind
of Siegel upper half space $\hat\Bbb H_2$ in (3) below, and the
other is the presentation of the ideas and explicit methods of
transforming the non-irreducible object $\hat\Bbb H_2$ to the
irreducible one $\Bbb H$; and the latter will have further
important applications in the forthcoming number theoretic
researches. By the way, we want to point out further that when the
action of a discrete group is considered there would arise some
new difficulties to be overcome.







Throughout this paper, we will use the following specific
notations. The letter $\ve$ is used to stand for $1$ or $-1$, that
is $\ve =\pm 1$. The letters $p$ and $q$ are always used to denote
the following matrices of order $2$,
$$
p:=\f1{\sqrt 2}\pmatrix 1 & -1 \\ 1 & 1 \endpmatrix, \quad
q:=\pmatrix 0 & 1 \\ 1 & 0 \endpmatrix .
$$
The capital letters $P$ and $Q$ are always used to denote the
block matrices
$$
P:=\pmatrix p & 0 \\ 0 & p \endpmatrix ,\quad Q:=\pmatrix q & 0 \\ 0 & q \endpmatrix.
$$
Then simple computations show that
$$
p^{-1}=\f{1}{\sqrt 2}\pmatrix 1 & 1 \\ -1 & 1 \endpmatrix, \quad P^{-1}=\pmatrix p^{-1} & 0 \\ 0 & p^{-1} \endpmatrix,
$$
and
$$q^2=I,\,\,\, Q^2=I.$$
Here and throughout this paper, $I$ is used to denote the identity
matrix of proper order $n\geq 1$, which is not necessarily the
same at different occurrences. We always use the capital letter
$J$ to denote the block matrix
$$
J:=\pmatrix 0 & I\\ -I & 0\endpmatrix .
$$
In this paper, $J$ is also assumed to be of order $4$, or
equivalently, the above blocks $I$ and $0$ are assumed to be of
order $2$. Then the well known symplectic group $Sp(4,\Bbb{R})$ of
order $4$, which will be denoted by $\Omega_2$ throughout this
paper, is as follows
$$
\Omega_2:=Sp(4,\Bbb{R})=\left\{M=\pmatrix A & B\\ C & D \endpmatrix \mid A,B,C,D \in \Bbb{R}^{(2,2)},\, \ ^{t}MJM=J\right\}.
$$
Notice that we clearly have
$$
P,\,Q\in \Omega_2.
$$
Recall that the Siegel upper half space $\Bbb{H}_{2}$ of order $2$
is defined as
$$
\Bbb{H}_2=\left\{Z\in\Bbb{C}^{(2,2)}\mid\  ^{t}Z=Z, \ Im
Z>0\right\}=\left\{Z=\pmatrix \tau_1 & z\\ z & \tau_2 \endpmatrix
\mid \tau_1,\tau_2,z\in \Bbb{C}, Im Z>0\right\}.
$$
And the action of $\Omega_2$ on $\Bbb{H}_2$ is defined by
$$
\aligned
f:\quad &\Omega_2\times \Bbb{H}_2 \rightarrow \Bbb{H}_2\\
&\,\,(M,\,Z)\,\mapsto W=f(M,\,Z)=M<Z>.
\endaligned\tag 1
$$
Here, and throughout this paper, we will always use the definition
$$
M<Z>:=(AZ+B)(CZ+D)^{-1},\tag 2
$$
for any $Z\in \Bbb{H}_2$ and $M=\pmatrix A & B\\
C & D
\endpmatrix\in \Omega_2$.

The materials of this paper are arranged as follows. In \S2 we
will first give the exact definition of our object we will work
with throughout our series, i.e., the definition of the new kind
of Siegel upper half space $\hat\Bbb H_2$ of order $2$. Then we
will give an initial formulation of the action group. \S3 is
arranged to give an alternative formulation of $\hat\Bbb H_2$ and
the corresponding action group. \S4 is devoted to the
investigation of the bi-holomorphic mappings of $\hat\Bbb H_2$.
The result together with the arguments in \S2 will give a complete
formulation of the action group. This then becomes the basis for
our further arguments. In \S5 we consider the reduced form of a
pair of points in $\hat\Bbb H_2$, which will be used to simplify
largely the formulation of our results in the following sections.
In \S6 the most important symplectic metric is built. This is a
cornerstone for the materials that follow. \S7 is devoted to
consider the geodesic line and the distance connecting two points.
This of course presents one of the most important intrinsic
feature of the so called geometry. In the last \S8 the
corresponding symplectic volume element is considered explicitly.

\head 2. The new kind of Siegel upper half space and the action
group
\endhead

Based on the well known Siegel upper half space $\Bbb{H}_{2}$ of
order $2$, we now give the definition of the most important object
$\hat{\Bbb{H}}_2$ in this paper:
$$
\hat{\Bbb{H}}_2:=\{ Z\in \Bbb{H}_2\mid Q<Z>=Z\}.\tag 3
$$
From now on, this $\hat{\Bbb{H}}_2$ will always be called the new
kind Siegel upper half space of order $2$, as expressed in the
title of this paper.
Note that, for any $Z=\pmatrix z_1 & z_2 \\
z_3 & z_4
\endpmatrix \in \Bbb{C}^{(2,2)}$, $Q<Z>=Z$ if and only if $qZ=Zq$,
i.e., $z_1=z_4,\,\,z_2=z_3$. Thus by the definition of
$\Bbb{H}_2$,
$$
\hat{\Bbb{H}}_2=\left\{Z=\pmatrix \tau & z\\ z & \tau \endpmatrix
\mid\  \tau,z\in\Bbb{C},\ {\text Im}\tau>|{\text Im} z|\right\}.
\tag 4
$$

\demo{Remark 1} From (4) one can see easily that the freedom of
the elements $Z$ in $\hat{\Bbb{H}}_2$ over $\Bbb{C}$ is $2$.
Further, in view of the form of the matrices $Z$, we call them
bi-symmetric. Recall that the well known Siegel upper half spaces
$\Bbb{H}_2$ and $\Bbb{H}_1$ have freedoms $3$ and $1$ over
$\Bbb{C}$ respectively. Thus from this freedom point of view, our
new kind of Siegel upper half space $\hat{\Bbb{H}}_2$ can be
viewed as an intermediate case between the cases of $\Bbb{H}_2$
and $\Bbb{H}_1$. Therefore, by comparing the classical outstanding
work of Siegel [S2] in 1936, it might be interesting to establish
a basis for the symplectic geometry of the new object
$\hat{\Bbb{H}}_2$.
\enddemo

The first main result in this paper is an explicit formulation of
the maximal subgroup of $\Omega_2$, which can act on
$\hat{\Bbb{H}}_2$ by group action.

\proclaim{Theorem 1}
Let $\hat{\Omega}_2$ be defined as
$$
\hat{\Omega}_2:=\{M\in \Omega_2\mid M<Z>\in\hat{\Bbb{H}}_2 \ \text{for all}\,\, Z\in\hat{\Bbb{H}}_2\}.
$$
Then we have
$$
\hat{\Omega}_2=\{M\in \Omega_2\mid MQ=\varepsilon QM\}.
$$
And so $\hat{\Omega}_2$ is the maximal subgroup of $\Omega_2$,
which can act on $\hat{\Bbb{H}}_2$ under the action given by
${(1)}$.
\endproclaim

To prove Theorem 1, we first give a preliminary lemma.

\proclaim{Lemma 1}
Let $M=\pmatrix A & B\\ C & D \endpmatrix\in \Omega_2$. Then
$$
M<Z>=Z\ \text{for all}\,\, Z\in\hat{\Bbb{H}}_2 \ \text{if and only
if}\  M=\ve I \ or\ M=\ve Q.
$$
\endproclaim

\demo{Proof} The sufficiency is obvious by definition. So we only
need to prove the necessity, i.e., we need to prove that $M=\ve  I
\ or\ =\ve  Q$ if $M<Z>=Z\ \text{for all}\,\,
Z\in\hat{\Bbb{H}}_2$. From $Z=M<Z>=(AZ+B)(CZ+D)^{-1}$ we get
$$
AZ+B=Z(CZ+D)=ZCZ+ZD.\tag 5
$$
Taking $Z=\tau I$ with $\tau \in \Bbb{C}$ and ${\text Im}\tau>0$,
which is clearly in $\hat{\Bbb{H}}_2$, then the above equality
becomes $\tau^2C+\tau(D-A)-B=0$, and this leads to
$$
B=C=0,\ D=A \tag 6
$$
by considering the limits as $\tau\to 0.$ Substituting these into
(5), we see that the matrix $A$ must satisfy
$$
AZ=ZA\tag 7
$$
for any $Z\in \hat{\Bbb{H}}_2$. Now, if one puts $A=\pmatrix
a_{11} & a_{12}\\ a_{21} & a_{22}\endpmatrix,\ Z=\pmatrix z &
\tau\\ \tau & z\endpmatrix$, then by $(7)$ and direct computation
we have $a_{11}\tau+a_{12}z=a_{11}\tau+a_{21}z$ and
$a_{11}z+a_{12}\tau=a_{12}\tau+a_{22}z$, which clearly imply that
$a_{12}=a_{21}$ and $a_{11}=a_{22}$ respectively. This then
enables us to assume that $A$ is of the form
$$A=\pmatrix a & b\\ b & a\endpmatrix$$
with $a,\,\,b\in {\Bbb{R}}$. Again, in view of $M\in \Omega_2$,
there holds $A\ ^{t}D-B\ ^{t}C=I$. So we also have $A\ ^{t}A=I.$
By this and direct computation, we easily obtain $a^2+b^2=1,\
ab=0$. So now there exist exactly two possibilities: one is
$b=0,\,a=\pm 1$ and we derive $A=\ve  I$, and the other is
$a=0,\,b=\pm 1$ and we derive $A=\ve  q$. This together with (6)
implies that $M=\ve  I$  or $ =\ve  Q$ as what we need. The proof
of lemma 1 is thus complete.
\enddemo

Now, we turn to the proof of Theorem 1. By the definition of
$\hat{\Omega}_2$ in Theorem 1, for any $M\in {\Omega}_2$, it is in
$\hat{\Omega}_2$ if and only if $ M<Z>\in\hat{\Bbb{H}}_2 \
\text{for all}\,\, Z\in\hat{\Bbb{H}}_2. $ But by the definition of
$\hat{\Bbb{H}}_2$, $M<Z>\in\hat{\Bbb{H}}_2$ if and only if
$qM<Z>=M<Z>q$, i.e., $qM<Z>q^{-1}=M<Z>$, or $Q<M<Z>>=M<Z>$ by
definition 1. Using the simple property of group action, one can
see easily that this last equality is also equivalent to
$\left(M^{-1}QM\right)<Z>=Z$. Thus by Lemma 1 we can derive
$M^{-1}QM=\ve  I$ or $=\ve  Q$. However, the first case of
$M^{-1}QM=\ve  I$ is impossible since this would imply $Q=1$ which
is clearly impossible. In other words, the set of $M\in
{\Omega}_2$ satisfying $M^{-1}QM=\ve  I$ is void. So there is no
contribution to $\hat{\Omega}_2$ from this kind of case. And thus
we can only have the latter case of $M^{-1}QM=\ve  Q$, i.e.,
$QM=\varepsilon MQ$, or in other words, the contribution to
$\hat{\Omega}_2$ of the $M\in {\Omega}_2$ comes exactly from the
latter case of $QM=\varepsilon MQ$. So this proves that $M\in
\hat{\Omega}_2$ iff $QM=\varepsilon MQ$, as desired by the first
part of Theorem 1. As for the other parts of the theorem, the
maximal property of $\hat{\Omega}_2$ is obvious from its
definition, and the remaining things can be derived easily from
the relative definitions by using the conclusion of the first
part. The proof of Theorem 1 is thus complete.

\head{3. An alternative formulation of $\hat\Bbb H_2$} and the
action group
\endhead

In this section, we first come to give an alternative formulation
of the new Siegel upper half space $\hat\Bbb{H}_2$, which is
isomorphic to $\hat\Bbb{H}_2$ under some "conformal
transformation", and consider the corresponding action group, see
$\hat E_2$ and $\hat\O_{\hat E_2}$ below in (10) and (12)
respectively. Then present another main result in this paper: The
action of $\O_2$ on $\hat\Bbb H_2$ is transitive. To this end, we
will first give some further conventions for notational
convenience. For any $m\times m$ matrix $A$ and any $m\times n$
matrix $X$, we use $\bar X$ to denote the conjugate of $X$, i.e.,
the matrix with all of its elements being the complex conjugates
of that of $X$, and we denote
$$
A\{X\}\!:=\! \ ^{t}XA\bar{X},\,\, A[X]\!:=\! \ ^{t}XAX.
$$
In this way we can write
$$
\align
\hat{\Bbb{H}}_2&=\left\{Z=\pmatrix \tau & z\\ z & \tau \endpmatrix \mid\  \tau,z\in\Bbb{C},\ {\text Im}\tau>|{\text Im} z|\right\}\\
&=\left\{Z\in\Bbb{C}^{(2,2)}\mid\ qZ=Zq,\ {\text Im}
Z=\frac{1}{2i}(Z-\bar{Z})>0\right\}.
\endalign
$$
And for any $Z\in\Bbb{C}^{(2,2)}$, in view of $J=\pmatrix 0 & I\\ -I & 0\endpmatrix$, we have by simple computation,
$$
J\bmatrix Z\\ I\endbmatrix := J\bmatrix\pmatrix Z\\ I\endpmatrix\endbmatrix
=(\ ^{t}Z\ I)J\pmatrix Z\\ I\endpmatrix=\ ^{t}Z-Z
$$
and
$$
J\left\{Z\atop{ I}\right\}:= J\left\{\pmatrix Z\\ I\endpmatrix\right\}
=(\ ^{t}Z\ I)J\pmatrix \bar{Z}\\ I\endpmatrix=\ ^{t}Z-\bar{Z}.
$$
These show that
$$
 ^{t}Z=Z \,\,\text{if and only if}\,\, J\bmatrix Z\\ I\endbmatrix=0,
$$
and when $J\bmatrix Z\\ I\endbmatrix=0$, we have
$$
\text{ Im} Z>0\,\,\text{if and only if}\,\,\frac{1}{2i}J\left\{Z\atop{ I}\right\}>0.
$$
Thus we can write
$$
\Bbb{H}_2=\left\{Z\in\Bbb{C}^{(2,2)}\mid\  J\bmatrix Z\\ I\endbmatrix=0,\ \frac{1}{2i}J\left\{Z\atop{ I}\right\}>0\right\},
$$
and
$$
\align
\hat{\Bbb{H}}_2 & =\left\{Z\in\Bbb{C}^{(2,2)}\mid\  J\bmatrix Z\\ I\endbmatrix=0,\ \frac{1}{2i}J\left\{Z\atop{ I}\right\}>0,\ qZ=Zq\right\}\\
& =\left\{Z\in\Bbb{C}^{(2,2)}\mid\ qZ=Zq ,\ \frac{1}{2i}J\left\{Z\atop{ I}\right\}>0\right\}.
\endalign
$$
Here the last equality comes from the fact that $qZ=Zq$ implies
$J\bmatrix Z\\ I\endbmatrix=0$. As for the action group
$\hat{\Omega}_2$ (sometimes also called motion group) of
$\hat{\Bbb{H}}_2$, we have
$$
\align
\hat{\Omega}_2 & =\left\{M\in \Omega_2 \mid \ QM=\pm MQ=\varepsilon MQ\right\}\\
&=\left\{ M\in \Bbb{R}^{(4,4)}\mid J[M]=J,\ J\{M\}=J,\ QM=\varepsilon MQ \right\}\\
&=\left\{ M\in \Bbb{R}^{(4,4)}\mid J[M]=J,\ QM=\varepsilon MQ
\right\}.\tag 8
\endalign
$$
Further, for $M=\pmatrix A & B\\ C & D\endpmatrix\in
\hat{\Omega}_2,\ Z\in\hat{\Bbb{H}}_2$, by putting $U:=AZ+B,\
V:=CZ+D$, then $M\pmatrix Z\\ I\endpmatrix=\pmatrix A & B\\ C &
D\endpmatrix\pmatrix Z\\ I\endpmatrix=\pmatrix U\\ V\endpmatrix$.
So
$$
J\bmatrix U\\ V\endbmatrix
=J\left[\pmatrix A & B\\ C & D\endpmatrix\pmatrix Z\\ I\endpmatrix\right]
=(Z\ I)\ ^{t}MJM\pmatrix Z\\ I\endpmatrix=J\bmatrix Z\\ I\endbmatrix=0,
$$
and
$$
\frac{1}{2i}J\left\{U\atop{ V}\right\}
=\frac{1}{2i}J\left\{\pmatrix A & B\\ C & D\endpmatrix\pmatrix Z\\ I\endpmatrix\right\}
=\frac{1}{2i}(Z\ I)\ ^{t}MJM\pmatrix \bar{Z}\\ I\endpmatrix
=\frac{1}{2i}J\left\{Z\atop{ I}\right\}>0.
$$
Thus
$$
^{t}UV-\ ^{t}VU=0,\ \frac{1}{2i}(^{t}U\bar{V}-\ ^{t}V\bar{U})>0.
$$
Now we come to establish a transformation which is similar to the so-called  "conformal transformation" of the plane of complex numbers. We claim at first that  $Z+iI$ is invertible for any $Z\in \hat{\Bbb{H}}_2$. In fact,
for any $v=v^{(2,1)}\in\Bbb{C}^{(2,1)}$ such that $(Z+iI)v=0$, one has
$-iv=Zv,\ i\bar{v}=\bar{Z}\bar{v},\ -i\ ^{t}v=\ ^{t}v\ ^{t}Z=\ ^{t}vZ$. So
$$
\frac{1}{2i}(Z-\bar{Z})\{v\}=\frac{1}{2i}\ ^{t}v(Z-\bar{Z})\bar{v}=\frac{1}{2i}(\ ^{t}vZ\bar{v}-\ ^{t}v\bar{Z}\bar{v})
=\frac{1}{2i}(-i\ ^{t}v\bar{v}-i\ ^{t}v\bar{v})=-\ ^{t}v\bar{v}\leq 0.
$$
On the other hand, since Im$Z$ is positive definite, we have
$$
\frac{1}{2i}(Z-\bar{Z})\{v\}=\ ^{t}v (\text{Im} Z)\bar{v}\geq 0.
$$
The combination of the above yields $\frac{1}{2i}(Z-\bar{Z})\{v\}=0,$ so gives rise to $v=0$. This then proves the invertibility of the matrix $Z+iI$ as being stated above. Now we can define the above mentioned "conformal transformation" of $\hat{\Bbb{H}}_2$ into $\Bbb{C}^{(2,2)}$ as follows
$$
\align
\psi :\quad\hat{\Bbb{H}}_2 & \longrightarrow  \Bbb{C}^{(2,2)}\\
Z & \longmapsto  Z_0: =\psi (Z):=(Z-iI)(Z+iI)^{-1}.\tag 9
\endalign
$$
From this definition it is not hard to see that $Z_0$ is bisymmetric and $I-Z_0\bar{Z}_0$ is positive definite and Hermitian.
The first thing is because $qZ_0=Z_0q$ which follows easily from $qZ=Zq$. The second is a consequence of the relevant arguments of Siegel [].
Thus we are naturally led to define a domain $\hat{E}_2$ in $\Bbb{C}^{(2,2)}$ as follows
$$
\hat{E}_2:=\left\{Z_0\in\Bbb{C}^{(2,2)}\mid\ qZ_0=Z_0q,\
I-Z_0\bar{Z}_0>0\right\}. \tag 10
$$
And (9) maps $\hat{\Bbb{H}}_2$ into $\hat{E}_2$. Conversely, if
$Z_0\in\hat{E}_2$, then Siegel [] has proved that $I-Z_0$ is
invertible. Thus we can also define the following map
$$
\align
\phi :\quad\hat{E}_2 & \longrightarrow  \Bbb{C}^{(2,2)}\\
Z_0 & \longmapsto Z:=\phi(Z_0):=i(I+Z_0)(I-Z_0)^{-1}. \tag 11
\endalign
$$
By this definition it is also easy to prove that
$Z\in\hat{\Bbb{H}}_2$ for any $Z_0\in \hat{E}_2$, thus (11) maps
$\hat{E}_2$ into $\hat{\Bbb{H}}_2$. Further, direct computations
show that the composition of (9) with (11) is the identity mapping
of $\hat{\Bbb{H}}_2$ and the composition of (11) with (9) is the
identity mapping of $\hat{E}_2$. Hence the mappings (9) and (11)
are all invertible and they are inverse mappings of each other,
and whence both of them are one-to-one correspondence. Therefore
the domain $\hat{E}_2$ defined by (10) can be served as another
formulation of our Siegel upper half space $\hat{\Bbb{H}}_2$. In
particular, one has
$$
\phi (0)=\psi ^{-1}(0)=iI,\ \psi (iI)=\phi ^{-1}(iI)=0.
$$

Next, we come to consider the action group on $\hat{E}_2$ corresponding to $\hat\Omega_2$. First of all, we define
$$
L:=\pmatrix iI & iI\\ -I & I\endpmatrix ,\ R:=\pmatrix -I & 0\\ 0 & I\endpmatrix .
$$
Note that there holds
$$
\align
&J[L]=^{t}LJL=2iJ,\ J[L^{-1}]=\frac{1}{2i}J,\ QL=LQ,\\
&\frac{1}{2i}J\{L\}=\frac{1}{2i}\ ^{t}LJ\bar{L}= \pmatrix -I & 0\\ 0 & I\endpmatrix = R.
\endalign
$$
Then for any $M\in \hat{\Omega}_2$, we define a corresponding matrix $M_0\in \Bbb{C}^{(4,4)}$ by
$$
M_0:=\pmatrix A_0 & B_0\\ C_0 & D_0\endpmatrix :=L^{-1}ML,
$$
where $A_0,\,B_0,\,C_0,\,D_0\in \Bbb{C}^{(2,2)}$, and we put
$$
\hat{\Omega}_{\hat{E}_2}
:=L^{-1}\hat{\Omega}_2L=\left\{M_0=L^{-1}ML \mid M\in
\hat{\Omega}_2\right\}. \tag 12
$$
Notice also that (11) is indeed the one to one correspondence
$$
\align
L:\ \hat{E}_2 & \longrightarrow \hat{\Bbb{H}}_2\\
Z_0 & \longmapsto Z=i(I+Z_0)(I-Z_0)^{-1}= L<Z_0>.\tag 13
\endalign
$$
This in combination with (1) implies that for any
$Z_0\in\hat{E}_2$ there holds
$$
M_0<Z_0>=(L^{-1}ML)<Z_0>=L^{-1}M<Z>=L^{-1}<W>.
$$
This shows that $W_0:=M_0<Z_0>$ is an element of $\hat{E}_2$, by
noting that $W=M<Z>$ is in $\hat{\Bbb{H}}_2$ since $Z$ is. Thus if
we note also that $\hat{\Omega}_{\hat{E}_2}$ is a group with
matrices multiplication, then it can be checked easily that we
have defined a group action of $\hat{\Omega}_{\hat{E}_2}$ on
$\hat{E}_2$ as follows
$$
\align
\hat{\Omega}_{\hat{E}_2} \times \hat{E}_2 & \longrightarrow \hat{E}_2 \\
\left(M_0,\,\,Z_0\right) & \longmapsto W_0=M_0<Z_0> . \tag 14
\endalign
$$
To have a better understanding of the group
$\hat{\Omega}_{\hat{E}_2}$ defined by (12), we need to give a more
explicit expression of it. To this end, we first note that by (8)
one can derive easily that
$$
\hat{\Omega}_2 =\left\{ M\in \Bbb{C}^{(4,4)}\mid J[M]=J,\ J\{M\}=J,\ QM=\varepsilon MQ \right\}.
$$
Thus for our purpose we only need to transform the constrains on
$M$ in this expression to that on $M_0$ in (12). This can be done
directly from $M=LM_0L^{-1}$. In deed, it is not difficult to find
by direct computations that $J[M]=J,\ J\{M\}=J$ and
$QM=\varepsilon MQ$ are equivalent to
$J[M_0]=J,\,R\left\{M_0\right\}=R$ and $QM_0=\varepsilon M_0Q$
respectively. Thus we have
$$
\hat{\Omega}_{\hat{E}_2}=L^{-1}\hat{\Omega}_2L=\{M_0\in\Bbb{C}^{(4,4)}\mid
J[M_0]=J,\ R\{M_0\}=R, \ QM_0=\varepsilon M_0Q\}. \tag 15
$$
Further, if we let
$$
F:=JR=\pmatrix 0 & I\\ I & 0\endpmatrix ,
$$
then for any invertible matrix
$
M_0=\pmatrix A_0 & B_0\\ C_0 & D_0\endpmatrix \in \Bbb{C}^{(4,4)}
$
with $A_0,\,B_0,\,C_0,\,D_0\in \Bbb{C}^{(2,2)}$ we can verify easily that
$$
M_0^{-1}F\bar{M_0}=F \,\text{ iff }\, F\bar{M_0}=M_0F \,\text{ iff }\, C_0=\bar{B}_0,\ D_0=\bar{A}_0.
$$
And, for this kind of $M_0$, we can also verify that $R\{M_0\}=R$ together with $J[M_0]=J$ implies
$F\bar{M_0}=M_0F$. This proves that
$$
\hat{\Omega}_{\hat{E}_2} \subseteq \left\{M_0=\left.\pmatrix A_0 & B_0\\ \bar{B}_0 & \bar{A}_0\endpmatrix \right| A_0,\,B_0 \in \Bbb{C}^{(2,2)},\,J[M_0]=J,\
 QM_0=\varepsilon M_0Q\right\}.
$$
Conversely, for any invertible matrix $M_0=\pmatrix A_0 & B_0\\ \bar{B_0} & \bar{A_0}\endpmatrix \in \Bbb{C}^{(4,4)}$ with $J[M_0]=J$, one can verify easily from $F\bar M_0=M_0F$ that $R\{M_0\}=R$. This proves that
$$
\left\{M_0=\left.\pmatrix A_0 & B_0\\ \bar{B}_0 & \bar{A}_0\endpmatrix \right| A_0,\,B_0 \in \Bbb{C}^{(2,2)},\,J[M_0]=J,\
 QM_0=\varepsilon M_0Q\right\} \subseteq \hat{\Omega}_{\hat{E}_2}.
$$
Gathering together the above we therefore obtain
$$
\aligned
\hat{\Omega}_{\hat{E}_2} & =\left\{M_0=\left.\pmatrix A_0 & B_0\\ \bar{B}_0 & \bar{A}_0\endpmatrix \right| A_0,\,B_0 \in \Bbb{C}^{(2,2)},\,J[M_0]=J,\ QM_0=\varepsilon M_0Q\right\}\\
& =\left\{M_0=\left.\pmatrix A_0 & B_0\\ \bar{B}_0 & \bar{A}_0\endpmatrix \right| A_0,\,B_0 \in \Bbb{C}^{(2,2)},\,
A_0\ ^{t}\bar{A}_0-B_0\ ^{t}\bar{B}_0=I,\ A_0\ ^{t}{B}_0=B_0\ ^{t}{A}_0 ,\ QM_0=\varepsilon M_0Q\right\}.
\endaligned \tag 16
$$

Next, we come to give a purely algebraic lemma which will be useful for our further arguments.

\proclaim{Lemma 2} Suppose that $K=\pmatrix k_1 & k_2 \\ k_2 & k_1 \endpmatrix$ is a positive definite real matrix of order $2$, then there exists a invertible real matrix $K_0=\pmatrix x_1 & x_2 \\ x_3 & x_4 \endpmatrix$ of order $2$ such that
$$K_0K \ ^{t}K_0=I,\quad qK_0=\varepsilon K_0q $$
with $\varepsilon =\pm 1 .$
\endproclaim

\demo{Proof} First of all, from the positive definiteness of $K$
we see that there holds $k_1>\left| k_2 \right|$. So we separate
the proof into two cases according to $k_2=0$ or not. If $k_2=0$,
the result is obvious by taking $K_0=\pm k_1^{-1/2}\pmatrix 1 & 0
\\ 0 & \varepsilon  \endpmatrix$ or $K_0=\pm k_1^{-1/2}\pmatrix 0
& 1 \\ \varepsilon & 0 \endpmatrix$. As for the case of $k_2\ne
0$, we note at first that the condition $qK_0=\varepsilon K_0q$ is
equivalent to $K_0$ being of the form $K_0=\pmatrix x_1 & x_2 \\
\varepsilon x_2 & \varepsilon x_1 \endpmatrix$ by direct
computation. Thus to prove the lemma, we are led to consider the
solvability of the matrix equation
$$
\pmatrix x_1 & x_2 \\ \varepsilon x_2 & \varepsilon x_1 \endpmatrix \pmatrix k_1 & k_2 \\ k_2 & k_1 \endpmatrix
\pmatrix x_1 & \varepsilon x_2 \\ x_2 & \varepsilon x_1 \endpmatrix =\pmatrix 1 & 0 \\ 0 & 1 \endpmatrix ,
$$
and again by direct computations, this can be shown to be
equivalent to the solvability of the system of the algebraic
equations
$$
\cases
k_1x_1^2+2k_2x_1x_2+k_1x_2^2=1 ,   \\
k_2x_1^2+2k_1x_1x_2+k_2x_2^2=0 .
\endcases
$$
However, this is clearly true since on noting $k_1>\left| k_2 \right|$ and $k_2\ne 0$ one can easily give the solutions of the system of the equations
as follows
$$
\cases
x_1=\frac{1}{2}\left(\varepsilon_1\left(k_1+k_2\right)^{-1/2}+\varepsilon_2\left(k_1-k_2\right)^{-1/2}\right)\\
x_2=\frac{1}{2}\left(\varepsilon_1\left(k_1+k_2\right)^{-1/2}-\varepsilon_2\left(k_1-k_2\right)^{-1/2}\right)
\endcases
$$
where $\varepsilon_1=\pm 1$, $\varepsilon_2=\pm 1$. The proof of Lemma 2 is complete.
\enddemo

Now we can state the main theorem in this section.

\proclaim{Theorem 2} The action of $\hat\Omega_{\hat E_2}$ on
$\hat E_2$ is transitive, so is the action of $\hat\Omega_2$ on
$\hat\Bbb{H}_2$.
\endproclaim

\demo{Proof} For the first assertion we only need to prove that
for any $Z_0\in\hat{E}_2$ it is in the same orbit of
$0\in\hat{E}_2$. Put $K=I-Z_0\bar{Z}_0$, which is clearly real,
positive definite, and bisymmetric. So by Lemma 2 there exists an
invertible real matrix $A_0$ such that
$A_0\left(I-Z_0\bar{Z}_0\right)\ ^{t}{A}_0=A_0K\ ^{t}{A}_0=I$ and
$qA_0=\varepsilon A_0q$. Thus if we let $B_0=-A_0Z_0,$ then by
(16) it is easy to verify that the matrix $M_0=\pmatrix A_0 &
B_0\\ \bar{B_0} & \bar{A_0}\endpmatrix$ is in
$\hat{\Omega}_{\hat{E}_2}$. Also, by the definition of $B_0$ we
clearly have $M_0<Z_0>=(A_0Z_0+B_0)(\bar B_0Z_0+\bar A_0)^{-1}=0$.
That is, for any given $Z_0\in\hat{E}_2$, there does exist $M_0\in
\hat{\Omega}_{\hat{E}_2}$ such that $M_0<Z_0>=0$ as desired by the
first assertion of our Theorem. To prove the second assertion, we
first take an arbitrary element $Z\in \hat\Bbb{H}_2$ and put
$Z_0=L^{-1}<Z>$. Then by the first assertion we can take an
$M_0\in \hat{\Omega}_{\hat{E}_2}$ such that $M_0<Z_0>=0$. Now by
taking $M=LM_0L^{-1}$, which is clearly in $\hat\Omega_2$, we can
obtain $M<Z>=iI$ by using the action of $L$ to both sides of
$M_0<Z_0>=0$ . This proves that $Z$ is in the orbit of $iI\in
\hat\Bbb H_2$ as desired. And thus the proof of Theorem 2 is
complete.
\enddemo

We now take a step further to consider the stability group of a point $Z_0$ in $\hat E_2$ under the action of $\hat\Omega_{\hat E_2}$, and then that of a point $Z$ in $\hat\Bbb H_2$ under the action of $\hat\Omega_2$. For any $Z_0 \in \hat E_2$, we use $S^{(1)}_{Z_0}$ to denote its stability group in $\hat\Omega_{\hat E_2}$, that is, we define
$$
S^{(1)}_{Z_0}:=\left\{\left.M_0\in \hat\Omega_{\hat E_2}\right|\
M_0<Z_0>=Z_0\right\}.\tag 17
$$
In particular, we have
$$
S^{(1)}_{0}:=\left\{\left.M_0\in \hat\Omega_{\hat E_2}\right|\
M_0<0>=0\right\}.\tag 18
$$
By Theorem 2 we know that there exists an $M_1\in \hat\Omega_{\hat
E_2}$ such that $M_1<0>=Z_0$, so the latter set in (18) does has
general meaning as the former set in (17).  More precisely, we
have the following proposition.

\proclaim{Proposition 1} Let $M_1$ be an element in $\hat\Omega_{\hat E_2}$ such that $M_1<0>=Z_0$. Then we have
$$
S^{(1)}_{Z_0}=M_1S^{(1)}_{0}M_1^{-1}=\left\{\left.M_1M_0M_1^{-1}\right|\
M_0\in S^{(1)}_{0}\right\}.\tag 19
$$
In particular, the right hand side of $(19)$ is irrelevant to the
choice of $M_1$ with $M_1<0>=Z_0$.
\endproclaim

\demo{Proof} This is just a direct consequence of the definition of the stability group $S^{(1)}_{Z_0}$ together with the definition of group action.
\enddemo

Based on Proposition 1, to understand more about the  stability group $S^{(1)}_{Z_0}$, one only needs to know more about $S^{(1)}_{0}$. As for this, we have the following

\proclaim{Proposition 2} Let $S^{(1)}_{0}$ be defined as in
$(18)$. Then we have
$$
\aligned
S^{(1)}_{0}& =\left\{\left.M_0=\pmatrix A_0 & 0 \\ 0 & \bar A_0 \endpmatrix\right|\ A_0\in \Bbb{C}^{(2,2)},\ A_0\ ^t\bar A_0=1,\ qA_0=\varepsilon A_0q\right\}\\
&=\left\{\left.M_0=\pmatrix A_0 & 0 \\ 0 & \bar A_0 \endpmatrix\right|\ A_0 = \pmatrix \left(\xi_1+\xi_2\right)/2 & \left(\xi_1-\xi_2\right)/2 \\
\varepsilon \left(\xi_1-\xi_2\right)/2 & \varepsilon \left(\xi_1+\xi_2\right)/2 \endpmatrix , \xi_1,\,\xi_2 \in \Bbb{C},\ \left|\xi_1\right|=\left|\xi_2\right|=1
\right\}.
\endaligned
$$
And thus the action of an element $M_0=\pmatrix A_0 & 0 \\ 0 & \bar A_0 \endpmatrix$ in $S^{(1)}_{0}$ on a point $Z_0$ in $\hat E_2$ becomes
$$
W_0:=M_0<Z_0>=A_0Z_0\bar A_0^{-1}=\ ^t UZ_0U,
$$
where $U=\ ^t A_0=\bar A_0^{-1}$ is a unitary matrix satisfying $qU=\varepsilon Uq$.
\endproclaim

\demo{Proof} The first equality for $S^{(1)}_{0}$ comes directly from the definition together with the observation that $QM_0=\varepsilon M_0Q$ is equivalent to $qA_0=\varepsilon A_0q$ if $M_0=\pmatrix A_0 & 0 \\ 0 & \bar A_0 \endpmatrix$. To prove the second equality, we first note that the condition $qA_0=\varepsilon A_0q$ on $A_0$ is equal to the assumption that $A_0$ is of the form $A_0=\pmatrix a & b \\ \varepsilon b & \varepsilon a \endpmatrix$ with $a,\,b \in \Bbb{C}$. For this kind of $A_0$, the condition $A_0\ ^t\bar A_0=1$ becomes both of the conditions $\left|a\right|^2+\left|b\right|^2=1$ and $a\bar b+b\bar a=0$, which is equivalent to $\left|a+b\right|^2=\left|a-b\right|^2=1$. Then denoting $\xi_1=a+b$ and $\xi_2=a-b$ we arrive at $a=\left(\xi_1+\xi_2\right)/2$ and $b=\left(\xi_1-\xi_2\right)/2$ with $\left|\xi_1\right|=\left|\xi_2\right|=1$ as desired. The proof of Proposition 2 is complete.
\enddemo

We now turn to the case of the action $\hat{\Omega}_2$ on
$\hat{\Bbb H}_2$. Similar to $(17)$, for any $Z\in \hat{\Bbb H}_2$
we define the stability group of $Z$ by
$$
S^{(2)}_{Z}:=\left\{\left.M\in \hat\Omega_{2}\right|\
M<Z>=Z\right\}.\tag 20
$$
And in particular
$$
S^{(2)}_{iI}:=\left\{\left.M\in \hat\Omega_{2}\right|\
M<iI>=iI\right\}.\tag 21
$$
Then by the mapping $(13)$ we can transform the above conclusions
to the following corresponding results without any difficulty.

\proclaim{Proposition 1$'$} For any given $Z\in \hat\Bbb{H}_2$, let $M_2$ be an element in $\hat\Omega_{2}$ such that $M_2<iI>=Z$. Then we have
$$
S^{(2)}_{Z}=M_2S^{(2)}_{iI}M_2^{-1}=\left\{\left.M_2MM_2^{-1}\right|\
M\in S^{(2)}_{iI}\right\}.\tag 22
$$
In particular, the right hand side of $(22)$ is irrelevant to the
choice of $M_2$ with $M_2<iI>=Z$.
\endproclaim

\proclaim{Proposition 2$'$} Let $S^{(2)}_{iI}$ be defined as in
$(21)$. Then we have
$$
S^{(2)}_{iI} = LS_0^{(1)}L^{-1} = \left\{LM_0L^{-1}\left|\ M_0\in S_0^{(1)}\right.\right\},
$$
where $S_0^{(1)}$ is as in Proposition $2$. And the action of an
element $M=LM_0L^{-1}$ in $S^{(2)}_{iI}$, with $M_0$ being of the
form $M_0=\pmatrix A_0 & 0 \\ 0 & \bar A_0 \endpmatrix$, on a
point $Z$ in $\hat\Bbb H_2$, is given by the formula
$$
L^{-1}<W>=A_0\left(L^{-1}<Z>\right)\bar A_0^{-1}=\ ^t U\left(L^{-1}<Z>\right)U,
$$
or identically,
$$
\left(W-iI\right)\left(W+iI\right)^{-1}=^t U\left(Z-iI\right)\left(Z+iI\right)^{-1}U,
$$
where $U=\ ^t A_0=\bar A_0^{-1}$ is again a unitary matrix satisfying $qU=\varepsilon Uq$.
\endproclaim

\demo{Remark 2} At the end of this section we remark that for each
$M_0\in \hat \Omega_{\hat E_2}$ it is not difficult to see the
mapping $W_0=M_0<Z_0>$ from $\hat E_2$ onto itself is
bi-holomorphic when $\hat E_2$ is considered as a domain of
$\Bbb{C}^{(1,2)}=\Bbb{C}\times\Bbb{C}$. In particularly, to each
$M_0\in S_0^{(1)}$ there corresponds to a bi-holomorphic mapping
$W_0=M_0<Z_0>$ from $\hat E_2$ onto itself with the fixed point
$0\in \hat E_2$. So $S_0^{(1)}$ can be viewed as a set of all
bi-holomorphic mappings $W_0=M_0<Z_0>$ with the fixed point $0\in
\hat E_2$ when $M_0\in \hat \Omega_{\hat E_2}$, and Proposition 2
gives a description of this kind of mappings. Accordingly, for
each $M\in \hat \Omega_{2}$ the mapping $W=M<Z>$ from $\hat\Bbb
H_2$ onto itself is bi-holomorphic when $\hat\Bbb H_2$ is
considered as a domain of $\Bbb{C}^{(1,2)}=\Bbb{C}\times\Bbb{C}$,
and thus Proposition 2$'$ can be viewed as a formulation of all
the bi-holomorphic mappings $W=M<Z>$ with the fixed point $iI\in
\hat\Bbb H_2$ when $M\in \hat \Omega_{2}$. Now, an important
converse question arises: what can we say about a general
bi-holomorphic mapping from $\hat E_2$ (resp. $\hat\Bbb H_2$) onto
itself, or more precisely, does every bi-holomorphic mapping from
$\hat E_2$ (resp. $\hat\Bbb H_2$) onto itself has the form of
$W_0=M_0<Z_0>$ (resp. $W=M<Z>$)? Similar to Siegel [S2], by the
relevant results in the preceding sections, this is equivalent to
asking specifically, wether or not are all the bi-holomorphic
mappings from $\hat E_2$ (resp. $\hat\Bbb H_2$) onto itself with
the fixed point $0\in \hat E_2$ (resp. $iI\in \hat\Bbb H_2$)
contained in $S_0^{(1)}$ (resp. $S_{iI}^{(2)}$)? These will be
answered in the following section.
\enddemo

\head 4. The group of bi-holomorphic mappings
\endhead

We first state a classical well known result from the analysis of several complex variables.

\proclaim{Lemma 3} Assume that $\left(r_1,\,r_2\right)$ is a permutation of $(1,\,2)$. Let $a_1,\,a_2\in \Bbb C$ satisfy $\left|a_1\right|<1,\,\left|a_2\right|<1$. Put $D:=\left\{z\in \Bbb C\ :\ \left|z\right|<1\right\}$. Then for any bi-holomorphic mapping
$$
\align
f:\ D\times D & \rightarrow D \times D \\
\left(z_1,\,z_2\right) & \mapsto \left(w_1,\,w_2\right)
\endalign
$$
satisfying $f\left(a_1,\,a_2\right)=f(0,\,0)$, it must be of the following form
$$
w_1=e^{i\theta_1}\frac{z_{r_1}-a_{r_1}}{1-\bar a_{r_1}z_{r_1}}, \quad w_2=e^{i\theta_2}\frac{z_{r_2}-a_{r_2}}{1-\bar a_{r_2}z_{r_2}},
$$
where $\theta_1,\,\theta_2$ are real parameters depending only on $f$ and satisfying $0\leq \theta_1,\,\theta_2<2\pi$.
\endproclaim

\demo{Proof} For a proof, one can see for example [H1].
\enddemo

The main result in this section is the following

\proclaim{Theorem 3} Every bi-holomorphic mapping from $\hat E_2$ (resp. $\hat\Bbb H_2$) onto itself has the form of $W_0=M_0<Z_0>$ (resp. $W=M<Z>$). So in view of Lemma 1, we can assert that the group of all bi-holomorphic mappings from $\hat E_2$ (resp. $\hat\Bbb H_2$) onto itself is exactly the quotient group $\hat\Omega_2/\{\pm I,\,\pm Q\}$ (resp. $\hat\Omega_{\hat E_2}/\{\pm I,\,\pm Q\})$.
\endproclaim

\demo{Proof} First of all, we note that, as pointed out at the end of last section, to prove the theorem, we only need to prove the following

Statement: every bi-holomorphic mapping $f$ from $\hat E_2$ onto itself with the fixed point $0\in \hat E_2$ is contained in $S_0^{(1)}$, that is, by Proposition 2, $f$ is of the form
$$
\align
f:\ \hat E_2&\longrightarrow \ \hat E_2\\
Z_0&\longmapsto\ f(Z_0)=\ ^{t}UZ_0U,
\endalign
$$
where $U$ is a unitary constant matrix satisfying $qU=\ve Uq$.

To prove the statement, we at first recall that by definition $p=\frac{1}{\sqrt{2}}\pmatrix 1 & -1\\ 1 & 1\endpmatrix $, and hence $^{t}p=\frac{1}{\sqrt{2}}\pmatrix 1 & 1\\ -1 & 1\endpmatrix=p^{-1}$. So for any $Z_0=\pmatrix z_1 & z_2\\ z_2 & z_1\endpmatrix\in\hat{E}_2$ we have
$$
^{t}pZ_0p=p^{-1}Z_0p=\pmatrix z_1+z_2 & 0\\ 0 & z_1-z_2\endpmatrix,
$$
and
$$
^{t}p(I-Z_0\bar{Z_0})p=I-p^{-1}Z_0pp^{-1}\bar{Z_0}p=\pmatrix 1-|z_1+z_2|^2 & 0\\ 0 & 1-|z_1-z_2|^2\endpmatrix.
$$
In particular, it is easy to see that the condition $I-Z_0\bar{Z_0}>0$, or $^{t}p(I-Z_0\bar{Z_0})p>0$ on $Z_0$, is equivalent to both of the conditions $1-|z_1+z_2|^2>0$ and $1-|z_1-z_2|^2>0$. Thus we can define a map $\sigma$ from $\hat{E}_2$ to $D\times D$ by
$$
\aligned
\sigma :\ \hat{E}_2&\longrightarrow\ D\times D\\
Z_0=\pmatrix z_1 & z_2\\ z_2 & z_1\endpmatrix &\longmapsto\ \sigma (Z_0):=(z_1+z_2, z_1-z_2)=(\hat{z}_1, \hat{z}_2),
\endaligned \tag 23
$$
where $\hat{z}_1=z_1+z_2$, $\hat{z}_2=z_1-z_2$. And it is not difficult to verify that $\sigma$ is bi-holomorphic and maps the zero $0=\pmatrix 0 & 0\\ 0 & 0 \endpmatrix \in \hat{E}_2$ to the zero $0=(0,\,0)\in D\times D$. Now define $\hat{f}$ to be $\hat{f}:=\sigma  f \sigma ^{-1}$. It is clear to find that $\hat{f}$ is a bi-holomorphic mapping from $D\times D$ onto itself with the fixed point $0=(0,0)$. Thus by Lemma 3 with $a_1=a_2=0$ we see that $\hat{f}$ must be one of the following two forms:
$$
\hat{f}(\hat{z}_1,\hat{z}_2)=(\xi_1\hat{z}_1,\ \xi_2\hat{z}_2),
$$
or
$$
\hat{f}(\hat{z}_1,\hat{z}_2)=(\xi_1\hat{z}_2,\ \xi_2\hat{z}_1),
$$
where $\left(\hat z_1,\,\hat z_2\right)$ denotes the variables in $D\times D$, and $\xi_1$ and $\xi_2$ are complex constant parameters depending only on $\hat{f}$ (so only on $f$) and satisfying $|\xi_1|=|\xi_2|=1$. Using $W_0=\pmatrix w_1 & w_2\\ w_2 & w_1\endpmatrix$ to denote $f(Z_0)$ and noting $\hat{f}=\sigma  f \sigma ^{-1}$, we can transform this to that of $f$ and obtain
$$
\align
&w_1=\frac{\xi_1\hat{z}_1+\xi_2\hat{z}_2}{2}=\frac{\xi_1+\xi_2}{2}z_1+\frac{\xi_1-\xi_2}{2}z_2,\\
&w_2=\frac{\xi_1\hat{z}_1-\xi_2\hat{z}_2}{2}=\frac{\xi_1-\xi_2}{2}z_1+\frac{\xi_1+\xi_2}{2}z_2
\endalign
$$
in the former case, and
$$
\align
&w_1=\frac{\xi_1+\xi_2}{2}z_1-\frac{\xi_1-\xi_2}{2}z_2,\\
&w_2=\frac{\xi_1-\xi_2}{2}z_1-\frac{\xi_1+\xi_2}{2}z_2
\endalign
$$
in the latter case. Both of them can clearly be unified to the form
$$
\aligned
&w_1=\frac{\xi_1+\xi_2}{2}z_1+\ve\frac{\xi_1-\xi_2}{2}z_2,\\
&w_2=\frac{\xi_1-\xi_2}{2}z_1+\ve\frac{\xi_1+\xi_2}{2}z_2
\endaligned \tag 24
$$
with $\ve =\pm 1$. Now we take $\eta_1=\xi_1^{1/2}$ and $\eta_2=\xi_2^{1/2}$ to be some fixed square-roots of $\xi_1$ and $\xi_2$ respectively, then put $u_1=\frac{\eta_1+\eta_2}{2},\ u_2=\frac{\eta_1-\eta_2}{2}$, and set
$$
U=\pmatrix u_1 & u_2\\ \ve u_2 & \ve u_1\endpmatrix.
$$
It is easy to see that $\eta_1$ and $\eta_2$ are complex constant parameters depending only on $f$ and satisfying $|\eta_1|=|\eta_2|=1$. So $U$ is a complex constant matrix depending only on $f$. And direct computations show that
$$
\align
\bar U\ ^t U & =\pmatrix \bar u_1 & \bar u_2 \\ \ve \bar u_2 & \ve \bar u_1 \endpmatrix \pmatrix u_1 & \ve u_2 \\ u_2 & \ve u_1 \endpmatrix
=\pmatrix \left| u_1\right|^2+\left| u_2\right|^2 & u_1 \bar u_2 +u_2 \bar u_1 \\ u_1 \bar u_2 +u_2 \bar u_1 & \left| u_1\right|^2+\left| u_2\right|^2 \endpmatrix \\
& = \pmatrix \left| \frac{\eta_1 + \eta_2}2 \right|^2+\left| \frac{\eta_1 - \eta_2}2 \right|^2 & \frac{\eta_1 + \eta_2}2 \frac{\bar\eta_1 - \bar\eta_2}2 +
\frac{\bar\eta_1 + \bar\eta_2}2\frac{\eta_1 - \eta_2}2 \\ \frac{\eta_1 + \eta_2}2 \frac{\bar\eta_1 - \bar\eta_2}2 +
\frac{\bar\eta_1 + \bar\eta_2}2\frac{\eta_1 - \eta_2}2 & \left| \frac{\eta_1 + \eta_2}2 \right|^2+\left| \frac{\eta_1 - \eta_2}2 \right|^2 \endpmatrix \\
& =\pmatrix 1 & 0 \\ 0 & 1 \endpmatrix = I,
\endalign
$$
which proves that $U$ is unitary. Further, in view of
$$
u_1^2+u_2^2=\left(\frac{\eta_1+\eta_2}2\right)^2+\left(\frac{\eta_1-\eta_2}2\right)^2=\frac{\xi_1+\xi_2}2
$$
and
$$
2u_1u_2=2\left(\frac{\eta_1+\eta_2}2\right)\left(\frac{\eta_1-\eta_2}2\right)=\frac{\xi_1-\xi_2}2,
$$
and then (24), it can be verified directly that
$$
\align
^t UZ_0U & = \pmatrix u_1 & \ve u_2 \\ u_2 & \ve u_1 \endpmatrix \pmatrix z_1 & z_2 \\ z_2 & z_1 \endpmatrix \pmatrix u_1 & u_2 \\ \ve u_2 & \ve u_1 \endpmatrix \\
& =  \pmatrix u_1z_1 + \ve u_2z_2 & u_1z_2 + \ve u_2z_1 \\ u_2z_1 + \ve u_1z_2 & \ve u_1z_1 + u_2z_2 \endpmatrix \pmatrix u_1 & u_2 \\ \ve u_2 & \ve u_1 \endpmatrix \\
& = \pmatrix \left(u_1^2 + u_2^2\right)z_1+2\ve u_1u_2z_2 & \ve \left(u_1^2 + u_2^2\right)z_2+2u_1u_2z_1 \\ \ve \left(u_1^2 + u_2^2\right)z_2+2u_1u_2z_1 & \left(u_1^2 + u_2^2\right)z_1+2\ve u_1u_2z_2 \endpmatrix \\
& = \pmatrix \frac{\xi_1+\xi_2}2z_1+\ve \frac{\xi_1-\xi_2}2z_2 & \frac{\xi_1-\xi_2}2z_1+\ve \frac{\xi_1+\xi_2}2z_2 \\
 \frac{\xi_1-\xi_2}2z_1+\ve \frac{\xi_1+\xi_2}2z_2 & \frac{\xi_1+\xi_2}2z_1+\ve \frac{\xi_1-\xi_2}2z_2 \endpmatrix \\
& = \pmatrix w_1 & w_2 \\ w_2 & w_1 \endpmatrix =W_0 =f(Z_0).
\endalign
$$
Again, by the definition of $U$ we can see plainly that $qU=\ve Uq$. Now, gathering together the above, we complete the proof of the above statement, and so complete the proof of Theorem 3.
\enddemo

\head 5. Reduced form for a pair of points of the space $\hat\Bbb
H_2$
\endhead

Recall that, in his remarkable paper [S2], based on the
transitivity of the action of $\Omega_2$ on the Siegel upper half
space $\Bbb H_2$ of order $2$, Siegel considered a simple reduced
form of any pair of points in the space $\Bbb H_2$ under the
action of the elements in $\Omega_2$. More precisely, he proved
that for any fixed pair of points $Z,\ Z_1\in \Bbb H_2$, there
exists an element $M\in \Omega_2$ such that both $M<Z_1>=iI$ and
$M<Z>=i\Lambda$ hold with $\Lambda$ being real and diagonal, and
satisfying $\Lambda \geq I$. Besides its independent interest from
the geometric point of view, this result then became a powerful
tool for the simplification of Siegel's relevant argument that
follows. In this section, based on our Theorem 2, similar to
Siegel[S2], we will take a step further to consider the reduced
form of any pair of points in the space $\hat \Bbb H_2$ under the
action of the elements in $\hat\Omega_2$. By philosophy and in
view of the transitivity of the action of $\hat\Omega_2$ on $\hat
\Bbb H_2$, this is roughly equal to finding a simple form of one
point in $\hat \Bbb H_2$ under the transformation of action by
element in $\hat\Omega_2$, which however must has the fixed point
$iI$. Our result is the following

\proclaim{Theorem 4} Suppose that $Z,\ Z_1$ is a fixed pair of
points in $\hat\Bbb H_2$, then there exists an element $M\in
\hat\Omega_2$ such that both $M<Z_1>=iI$ and $M<Z>=i\Lambda$ hold
with $\Lambda=\pmatrix \lambda_1 & \lambda_2 \\ \lambda_2 &
\lambda_1 \endpmatrix$ being real and bi-symmetric, and satisfying
$\lambda_1\geq \lambda_2 +1,\ \lambda_2 \geq 0$. Moreover,
$\Lambda$ is unique in the sense that if there exists another
element $M'\in \hat\Omega_2$ such that both $M'<Z_1>=iI$ and
$M'<Z>=i\Lambda '$ hold with $\Lambda '$ having the same
properties as that of $\Lambda$, then we have $\Lambda '=\Lambda$.
\endproclaim

\demo{Proof} By the transitivity of the action of $\hat \Omega_2$
on $\hat\Bbb H_2$, we can assume without loss of generality that
$Z_1=iI$. To see this, we only need to show that Theorem 4 is true
if one assumes its validity when $Z_1=iI$. In fact, for any given
$Z_1$ and $Z$, by transitivity it is known that there exists $M_1
\in \hat\Omega_2$ such that $M_1<Z_1>=iI$. Put $M_1<Z>=Z'$. Then
Theorem 4 for the pair of points $iI$ and $Z'$ tells us that there
exists $M_0 \in \hat\Omega_2$ such that $M_0<iI>=iI$ and
$M_0<Z'>=i\Lambda$. Taking $M=M_0M_1$, we can see easily that
$M<Z_1>=iI$ and $M<Z>=i\Lambda$. This is the existence part of
Theorem 4 in general. To see the uniqueness, we assume there exist
$M,\ M' \in \hat\Omega_2$ such that
$$
M<Z_1>=iI,\ M<Z>=i\Lambda
$$
and
$$
M'<Z_1>=iI,\ M'<Z>=i\Lambda '.
$$
These together with $M_1<Z_1>=iI$ and $M_1<Z>=Z'$ lead to
$$
MM_1^{-1}<iI>=iI,\ MM_1^{-1}<Z'>=i\Lambda
$$
and
$$
M'M_1^{-1}<iI>=iI,\ M'M_1^{-1}<Z'>=i\Lambda '.
$$
And this implies $\Lambda =\Lambda '$ by the uniqueness part of
Theorem 4 for the pair of numbers $iI$ and $Z'=M_1<Z>$, which is
what we need. Now we come to prove Theorem 4 under the assumption
that $Z_1=iI$. First, consider the existence part. Recalling (13)
and (23), we see that there exists a bijective mapping $\sigma
L^{-1}$ from $\hat \Bbb H_2$ to $D\times D$ which maps $iI$ to
$(0,\ 0)$ in particular. Let $\left(\hat z_1^0,\ \hat
z_2^0\right)$ denote the image of $Z$ under this mapping, i.e.,
$\left(\hat z_1^0,\ \hat z_2^0\right)=\left(\sigma
L^{-1}\right)(Z)$, and let $r_1=\max \left\{\left|\hat
z_1^0\right|,\ \left|\hat z_2^0\right|\right\}$, $r_2=\min
\left\{\left|\hat z_1^0\right|,\ \left|\hat z_2^0\right|\right\}$.
Then in view of Lemma 3 we can take a bi-holomorphic function $f$
with fixed point $(0,\ 0)$ from $D\times D$ to itself such that
$$
f\left(\hat z_1^0,\ \hat z_2^0\right)=\left(r_1,\ r_2\right). \tag
25
$$
Indeed, if we put $\hat z_j^0=\left|\hat z_j^0\right|{\text e}^{i\theta_j}$ for $1\leq j\leq 2$ with  $0\leq \theta_j <2\pi$, then when $\left|\hat z_1^0\right| \geq \left|\hat z_2^0\right|$, $f$ may be taken as
$$
f\left(\hat z_1,\ \hat z_2\right)=\left({\text e}^{-i\theta_1}\hat z_1,\ {\text e}^{-i\theta_2}\hat z_2\right),
$$
and when $\left|\hat z_1^0\right| < \left|\hat z_2^0\right|$, $f$ may be taken as
$$
f\left(\hat z_1,\ \hat z_2\right)=\left({\text e}^{-i\theta_2}\hat z_2,\ {\text e}^{-i\theta_1}\hat z_1\right),
$$
where $\left(\hat z_1,\ \hat z_2\right)$ is an arbitrary point in
$D\times D$. For this $f$, there corresponds to a bi-holomorphic
mapping $g=L\sigma ^{-1}f\sigma L^{-1}$ from $\hat \Bbb H_2$ to
itself, which clearly fixes the point $iI$. And since $f$ moves
the point $\left(\hat z_1^0,\ \hat z_2^0\right)$ to $\left(r_1,\
r_2\right)$ in $D\times D$ by (25), we can see easily that $g$
moves the corresponding point $Z=\left(\sigma
L^{-1}\right)^{-1}\left(\hat z_1^0,\ \hat z_2^0\right)$ of
$\left(\hat z_1^0,\ \hat z_2^0\right)$ to the corresponding point
$\left(\sigma L^{-1}\right)^{-1}\left(r_1,\ r_2\right)$ of
$\left(r_1,\ r_2\right)$ under the mapping $\sigma L^{-1}$. Then
we can take $M$ to be an element in $\hat \Omega _2$, for which
the action on $\hat \Bbb H_2$ is exactly identical to the mapping
$g$. Thus $M$ fixes $iI$, and by (23) and (13),
$$
\align
M<Z>&=g(Z)=\left(\sigma L^{-1}\right)^{-1}\left(r_1,\ r_2\right)=L\sigma^{-1}\left(r_1,\ r_2\right)=L\left(\pmatrix \f{r_1+r_2}2 &\f{r_1-r_2}2\\ \f{r_1-r_2}2 & \f{r_1+r_2}2\endpmatrix\right)\\
&=i\left(I+\pmatrix \f{r_1+r_2}2 &\f{r_1-r_2}2\\ \f{r_1-r_2}2 & \f{r_1+r_2}2\endpmatrix\right)
\left(I-\pmatrix \f{r_1+r_2}2 &\f{r_1-r_2}2\\ \f{r_1-r_2}2 & \f{r_1+r_2}2\endpmatrix\right)^{-1}
=i\Lambda ,
\endalign
$$
where
$$
\align
\Lambda &=\left(I+\pmatrix \f{r_1+r_2}2 &\f{r_1-r_2}2\\ \f{r_1-r_2}2 & \f{r_1+r_2}2\endpmatrix\right)
\left(I-\pmatrix \f{r_1+r_2}2 &\f{r_1-r_2}2\\ \f{r_1-r_2}2 & \f{r_1+r_2}2\endpmatrix\right)^{-1}\\
&=\pmatrix \f{1-r_1r_2}{\left(1-r_1\right)\left(1-r_2\right)} & \f{r_1-r_2}{\left(1-r_1\right)\left(1-r_2\right)} \\
\f{r_1-r_2}{\left(1-r_1\right)\left(1-r_2\right)} &
\f{1-r_1r_2}{\left(1-r_1\right)\left(1-r_2\right)}
\endpmatrix ,
\endalign
$$
which clearly satisfies the desired conditions for $\L$ in Theorem
4, and the existence part is proved. Next, turn to the proof of
the uniqueness of $\L$. So we assume there is another $M'$ in
$\hat \Omega_2$ such that $M'<iI>=iI$ and $M'<Z>=i\L '$, and we
are going to prove $\L =\L '$. From this assumption it can be seen
easily that the action of the element $M'M^{-1}\in \hat \Omega_2$
fixes the point $iI$ and satisfies $\left(M'M^{-1}\right)<i\L
>=i\L '$. Thus by (13) we get a bi-holomorphic mapping
$L^{-1}\left(M'M^{-1}\right)L$ from $\hat E_2$ onto itself such
that
$$
\left(L^{-1}\left(M'M^{-1}\right)L\right)\left<L^{-1}(i\L )\right>=L^{-1}(i\L '),
$$
which follows from $\left(M'M^{-1}\right)<i\L >=i\L '$. Hence by
Theorem 3 we can get a  unitary $2\times 2$ matrix $U$ with
$qU=\varepsilon Uq$ such that
$$
^t UL^{-1}(i\L )U=L^{-1}(i\L '). \tag 26
$$
Again by (13), we have by simple computations,
$$
L^{-1}(i\L )=\pmatrix
\f{\lambda_1^2-\lambda_2^2-1}{(\lambda_1+1)^2-\lambda_2^2} &
\f{2\lambda_2}{(\lambda_1+1)^2-\lambda_2^2} \\
\f{2\lambda_2}{(\lambda_1+1)^2-\lambda_2^2} &
\f{\lambda_1^2-\lambda_2^2-1}{(\lambda_1+1)^2-\lambda_2^2}\endpmatrix
=\pmatrix \f{r_1+r_2}2 & \f{r_1-r_2}2 \\ \f{r_1-r_2}2 &
\f{r_1+r_2}2 \endpmatrix ,
$$
and
$$
L^{-1}(i\L ')=\pmatrix \f{{\lambda '}_1^2-{\lambda
'}_2^2-1}{({\lambda '}_1+1)^2-{\lambda '}_2^2} &
\f{2\lambda'_2}{(\lambda_1+1)^2-\lambda_2^2} \\
\f{2\lambda'_2}{(\lambda_1+1)^2-\lambda_2^2} & \f{{\lambda
'}_1^2-{\lambda '}_2^2-1}{({\lambda '}_1+1)^2-{\lambda
'}_2^2}\endpmatrix =\pmatrix \f{r'_1+r'_2}2 & \f{r'_1-r'_2}2 \\
\f{r'_1-r'_2}2 & \f{r'_1+r'_2}2 \endpmatrix ,
$$
where we have used the symbols
$$
\Lambda=\pmatrix \lambda_1 & \lambda_2\\ \lambda_2 & \lambda_1\endpmatrix , \,\, \Lambda '=\pmatrix {\lambda '}_1 & {\lambda '}_2\\ {\lambda '}_2 & {\lambda '}_1\endpmatrix ,
$$
and
$$
r_1=\frac{\lambda_1+\lambda_2-1}{\lambda_1+\lambda_2+1},\ r_2=\frac{\lambda_1-\lambda_2-1}{\lambda_1-\lambda_2+1},\
r'_1=\frac{{\lambda '}_1+{\lambda '}_2-1}{{\lambda '}_1+{\lambda '}_2+1},\ r'_2=\frac{{\lambda '}_1-{\lambda '}_2-1}{{\lambda '}_1-{\lambda '}_2+1}.
$$
Under these notations, (26) can thus be rewritten as
$$
^t U \pmatrix \f{r_1+r_2}2 & \f{r_1-r_2}2 \\ \f{r_1-r_2}2 &
\f{r_1+r_2}2 \endpmatrix U= \pmatrix \f{r'_1+r'_2}2 &
\f{r'_1-r'_2}2 \\ \f{r'_1-r'_2}2 & \f{r'_1+r'_2}2 \endpmatrix ,
\tag 27
$$
and by the assumptions $\lambda_1\geq \lambda_2+1\geq 1$ and ${\lambda '}_1\geq {\lambda '}_2+1\geq 1$, we also have
$$
0\leq r_2\leq r_1<1, \  0\leq r'_2\leq r'_1<1,
$$
and
$$
\lambda_1=\frac{1-r_1r_2}{(1-r_1)(1-r_2)},\ \lambda_2=\frac{r_1-r_2}{(1-r_1)(1-r_2)},\ \lambda '_1=\frac{1-r'_1r'_2}{(1-r'_1)(1-r'_2)},\ \lambda '_2=\frac{r'_1-r'_2}{(1-r'_1)(1-r'_2)}.
$$
So, to prove $\L = \L '$, it is sufficient to prove $r_1=r'_1$ and
$r_2=r'_2$ from (27). We separate two cases according to $\ve =1$
or $=-1$ to do this. If $\ve=1$, then by $qU=Uq$ we see that
$U=\pmatrix u_1 & u_2 \\ u_2 & u_1\endpmatrix$ with $u_1,\,\,u_2
\in \Bbb C$. And hence by $\ ^{t}U\bar{U}=I$, we get
$$
|u_1|^2+|u_2|^2=1,\ u_1\bar{u}_2+u_2\bar{u}_1=0.\tag 28
$$
Again, multiplying on both sides of (27) by $^t p=p^{-1}$ and $p$
from the left and right respectively, and inserting $pp^{-1}$ into
the left-hand side of it behind $^t U$ and preceding $U$, we get
$$
\ ^{t}p\ ^{t}Upp^{-1}\pmatrix\frac{r_1+r_2}{2} & \frac{r_1-r_2}{2}\\ \frac{r_1-r_2}{2} & \frac{r_1+r_2}{2}\endpmatrix p\ ^{t}pUp
=p^{-1}\pmatrix\frac{r'_1+r'_2}{2} & \frac{r'_1-r'_2}{2}\\ \frac{r'_1-r'_2}{2} & \frac{r'_1+r'_2}{2}\endpmatrix p ,
$$
which gives
$$
\pmatrix u_1+u_2 & 0 \\ 0 & u_1-u_2\endpmatrix \pmatrix r_1 & 0 \\ 0 & r_2\endpmatrix\pmatrix u_1+u_2 & 0 \\ 0 & u_1-u_2\endpmatrix
=\pmatrix r'_1 & 0 \\ 0 & r'_2\endpmatrix .
$$
Whence we get
$$
r_1(u_1+u_2)^2=r'_1,\ r_2(u_1-u_2)^2=r'_2.
$$
Taking absolute values on both sides of these two equalities, we
then obtain by (28),
$$
r'_1=r_1\left|u_1+u_2\right|^2=r_1\left(\left|u_1\right|^2+\left|u_2\right|^2+u_1\bar u_2+u_2\bar u_1\right)=r_1,
$$
and similarly
$$
r'_2=r_2\left|u_1-u_2\right|^2=r_2\left(\left|u_1\right|^2+\left|u_2\right|^2-u_1\bar u_2-u_2\bar u_1\right)=r_2,
$$
as desired. If $\ve=-1$, then by $qU=-Uq$ we see that $U=\pmatrix
u_1 & u_2 \\ -u_2 & -u_1\endpmatrix$ again with $u_1$ and $u_2$
being complex numbers. Also using $\ ^{t}U\bar{U}=I$ we can get
(28). Then again by (27) we can obtain by direct computations
$$
r'_1=r_2(u_1+u_2)^2,\ r'_2=r_1(u_1-u_2)^2.
$$
Then taking absolute values and using (28), we obtain $r'_1=r_2,\
r'_2=r_1$. This together with the conditions $0\leq r_2\leq r_1<1$
and $0\leq r'_2\leq r'_1<1$ gives rise to
$$
r'_1\geq r'_2=r_1\geq r_2=r'_1,
$$
which clearly implies $r'_1=r_1=r'_2=r_2$ again as desired. The proof of Theorem 4 is complete.
\enddemo

\head 6. The symplectic metric for $\hat\Bbb H_2$
\endhead

The main purpose of this section is to consider the existence and
uniqueness, up to a positive constant multiple, of symplectic
metric for the space $\hat\Bbb H_2$ invariant under the action of
the group $\hat\Omega_2$. To this end, we first give a definition
of the cross ratio of points in $\hat\Bbb H_2$ following
Siegel[S2]. For any $Z,\,Z_1\in\hat{\Bbb H}_2$, in view of Im$Z>0$
and Im$Z_1>0$, we can see easily that the matrix $(Z-\bar{Z}_1)$
is invertible. So we can define the matrix
$$
{\frak
R}(Z,Z_1):=(Z-Z_1)(Z-\bar{Z}_1)^{-1}(\bar{Z}-\bar{Z}_1)(\bar{Z}-Z_1)^{-1}
, \tag 29
$$
and call it the cross ratio of $Z$ and $Z_1$. For any $M=\pmatrix A & B \\ C & D \endpmatrix \in\hat{\Omega}_2$, write
$$
W=M<Z>,\ W_1=M<Z_1>.
$$
Then direct computation yields
$$
\align
Z_1-Z&=\pmatrix -I & \!\!\! Z_1\endpmatrix\pmatrix Z \\ I\endpmatrix
=\pmatrix Z_1 & \!\!\! I\endpmatrix\pmatrix 0 & I\\ -I & 0\endpmatrix\pmatrix Z \\ I\endpmatrix
=\pmatrix Z_1 & \!\!\! I\endpmatrix\ ^{t}MJM \pmatrix Z \\ I\endpmatrix\\
&=\left(\pmatrix Z_1 & \!\!\! I\endpmatrix \pmatrix ^t A & ^t C \\ ^t B & ^t D \endpmatrix \right)
\pmatrix 0 & I\\ -I & 0\endpmatrix \left( \pmatrix A & B \\ C & D \endpmatrix \pmatrix Z \\ I\endpmatrix \right) \\
& = \pmatrix Z_1\, ^t\!A+ ^t\!B, & \!\!\! Z_1\, ^t\!C+ ^t\!D \endpmatrix \pmatrix 0 & I\\ -I & 0\endpmatrix \pmatrix AZ+B \\ CZ+D \endpmatrix \\
& = \pmatrix -\left(Z_1\, ^t\!C+ ^t\!D\right), & \!\!\! Z_1\, ^t\!A+ ^t\!B \endpmatrix \pmatrix AZ+B \\ CZ+D \endpmatrix \\
& = \ ^t\!\left(AZ_1+B\right)(CZ+D)- ^t\!\left(CZ_1+D\right)(AZ+B)\\
& = \ ^{t}(CZ_1+D)(W_1-W)(CZ+D).
\endalign
$$
Here for the last equality we have used the symmetric property of
$^{t}(CZ_1+D)(AZ_1+B)$ which follows from $q\,
^t(CZ_1+D)(AZ_1+B)=\, ^t(CZ_1+D)(AZ_1+B)q$. Using $\bar{Z}_1$ to
replace $Z_1$ in the above equality, we have
$$
\bar{Z}_1-Z=\ ^{t}(C\bar{Z}_1+D)(\bar{W}_1-W)(CZ+D).
$$
Substituting these two expressions into (29), we get
$$
\frak R(Z,Z_1)=\ ^{t}(CZ_1+D)\frak R(W,W_1)\ ^{t}(CZ_1+D)^{-1}.
\tag 30
$$
As a consequence, this shows that $\frak R\left(Z,Z_1\right)$ and
$\frak R\left(W,W_1\right)$ have the same eigenvalues. Conversely,
for two given pairs $\left(Z,\, Z_1\right)$ and $\left(W,\,
W_1\right)$, by Theorem 4 we realize that there exist elements
$M_1,\,M_2\in \hat\Omega_2$ such that $M_1<Z_1>=iI$,
$M_1<Z>=i\Lambda_1$, and $M_2<W_1>=iI$, $M_2<W>=i\Lambda_2$, with
$\Lambda_1=\pmatrix \lambda_1 & \lambda_2
\\ \lambda_2 & \lambda_1 \endpmatrix$ and $\Lambda_2=\pmatrix \lambda '_1 & \lambda '_2
\\ \lambda '_2 & \lambda '_1 \endpmatrix$, where
 $\lambda_1\geq \lambda_2 +1,\ \lambda_2 \geq 0$
and $\lambda '_1\geq \lambda '_2 +1,\ \lambda '_2 \geq 0$. Thus,
if $\frak R\left(Z,Z_1\right)$ and $\frak R\left(W,W_1\right)$
have the same eigenvalues, then $\frak R\left(iI,i\L_1\right)$ and
$\frak R\left(iI,i\L_2\right)$ have the same eigenvalues. So in
view of
$$
p^{-1}\frak R\left(Z,Z_1\right)p=\pmatrix
\left(\f{\l_1+\l_2-1}{\l_1+\l_2+1}\right)^2 & 0 \\ 0 &
\left(\f{\l_1-\l_2-1}{\l_1-\l_2+1}\right)^2 \endpmatrix
$$
and
$$
p^{-1}\frak R\left(W,W_1\right)p=\pmatrix
\left(\f{\l'_1+\l'_2-1}{\l'_1+\l'_2+1}\right)^2 & 0 \\ 0 &
\left(\f{\l'_1-\l'_2-1}{\l'_1-\l'_2+1}\right)^2 \endpmatrix
$$
we can realize
$$
\l_1=\l'_1,\,\,\l_2=\l'_2.
$$
Whence $\L_1=\L_2,$ and hence we have
$$
\left(M_2^{-1}M_1\right)<Z_1>=W_1,\,\,\left(M_2^{-1}M_1\right)<Z>=W.
$$
In summary, we have the following

\proclaim{Theorem 5} There exists a symplectic transformation
$M\in \hat\O_2$ mapping a given pair $\left(Z,\,Z_1\right)$ of
$\hat\Bbb{H}_2$ into another given pair $\left(W,\,W_1\right)$ of
$\hat\Bbb{H}_2$, if and only if the cross ratios $\frak R
\left(Z,\,Z_1\right)$ and $\frak R \left(W,\,W_1\right)$ have the
same eigenvalues.
\endproclaim

Again, (30) shows that the trace $tr(\frak R(Z,Z_1))$ of the cross
ratio $\frak R(Z,Z_1)$ is invariant under the action by $M\in
\hat{\Omega}_2$, i.e., we have for any $M\in\hat{\Omega}_2$,
$$
tr(\frak R(W,W_1))=tr(\frak R(M<Z>,M<Z_1>))=tr(\frak R(Z,Z_1)).
\tag 31
$$
Now we fix $Z$ and $M$ temporarily, and regard the right  and the
left hand sides of (31) as functions of $Z_1$ and $W_1$
respectively. Then by definition (29) we can find easily that
$\frak R(Z,Z_1)$ and so $tr(\frak R(Z,Z_1))$ are differentiable
with respect to the variable $Z_1$. Thus by (31) and in view of
the linearity of the differential operator, we get
$$
tr(d^2\frak R(W,W_1))=tr(d^2\frak R(Z,Z_1)). \tag 32
$$
In particular, by taking $Z_1=Z$, (32) shows that the trace of the
second order differential $d^2\frak R(Z,Z_1)$ at $Z$, which is
denoted by $\left.d^2\frak R(Z,Z_1)\right|_{Z_1=Z}$ as usual, is
invariant under the action by any $M\in\hat{\Omega}_2$. Again,
from (29), direct computation gives
$$
\align
{\text d}{\frak R}(Z,Z_1) = & {\text d}\left( (Z-Z_1)(Z-\bar{Z}_1)^{-1}(\bar{Z}-\bar{Z}_1)(\bar{Z}-Z_1)^{-1} \right)\\
= & -\left({\text d}Z_1\right) (Z-\bar{Z}_1)^{-1}(\bar{Z}-\bar{Z}_1)(\bar{Z}-Z_1)^{-1} \\
& + (Z-Z_1)(Z-\bar{Z}_1)^{-1} \left({\text d}\bar{Z}_1\right)(Z-\bar{Z}_1)^{-1} (\bar{Z}-\bar{Z}_1)(\bar{Z}-Z_1)^{-1}\\
& - (Z-Z_1)(Z-\bar{Z}_1)^{-1}\left({\text d}\bar{Z}_1\right)(\bar{Z}-Z_1)^{-1} \\
& + (Z-Z_1)(Z-\bar{Z}_1)^{-1}(\bar{Z}-\bar{Z}_1)(\bar{Z}-Z_1)^{-1}
\left({\text d}Z_1\right)(\bar{Z}-Z_1)^{-1} ,
\endalign
$$
and hence
$$
\align
{\text d}^2{\frak R}(Z,Z_1)= & -\left({\text d}Z_1\right)\left\{ (Z-\bar{Z}_1)^{-1} \left({\text d}\bar Z_1\right)
(Z-\bar{Z}_1)^{-1}(\bar{Z}-\bar{Z}_1)(\bar{Z}-Z_1)^{-1}\right.\\
& \quad - (Z-\bar{Z}_1)^{-1}\left({\text d}\bar Z_1\right)(\bar{Z}-Z_1)^{-1}\\
& \quad + \left. (Z-\bar{Z}_1)^{-1}(\bar{Z}-\bar{Z}_1)(\bar{Z}-Z_1)^{-1}\left({\text d}Z_1\right)(\bar{Z}-Z_1)^{-1} \right\}\\
& + \left\{ - \left({\text d}Z_1\right)(Z-\bar{Z}_1)^{-1} \left({\text d}\bar{Z}_1\right)(Z-\bar{Z}_1)^{-1} (\bar{Z}-\bar{Z}_1)(\bar{Z}-Z_1)^{-1} \right.\\
& \quad +(Z-Z_1)(Z-\bar{Z}_1)^{-1}\left({\text d}\bar{Z}_1\right)(Z-\bar{Z}_1)^{-1}\left({\text d}\bar{Z}_1\right)(Z-\bar{Z}_1)^{-1} (\bar{Z}-\bar{Z}_1)(\bar{Z}-Z_1)^{-1} \\
& \quad +(Z-Z_1)(Z-\bar{Z}_1)^{-1} \left({\text d}\bar{Z}_1\right)(Z-\bar{Z}_1)^{-1}\left({\text d}\bar{Z}_1\right)(Z-\bar{Z}_1)^{-1}(\bar{Z}-\bar{Z}_1)(\bar{Z}-Z_1)^{-1} \\
& \quad +(Z-Z_1)(Z-\bar{Z}_1)^{-1} \left({\text d}\bar{Z}_1\right)(Z-\bar{Z}_1)^{-1} (-{\text d}\bar{Z}_1)(\bar{Z}-Z_1)^{-1}\\
& \quad +\left.(Z-Z_1)(Z-\bar{Z}_1)^{-1} \left({\text d}\bar{Z}_1\right)(Z-\bar{Z}_1)^{-1} (\bar{Z}-\bar{Z}_1)(\bar{Z}-Z_1)^{-1}({\text d}Z_1)(\bar{Z}-Z_1)^{-1}\right\}\\
& - \left\{- \left({\text d}Z_1\right)(Z-\bar{Z}_1)^{-1}\left({\text d}\bar{Z}_1\right)(\bar{Z}-Z_1)^{-1} \right.\\
& \quad + (Z-Z_1)(Z-\bar{Z}_1)^{-1} \left({\text d}\bar Z_1\right) (Z-\bar{Z}_1)^{-1} \left({\text d}\bar{Z}_1\right)(\bar{Z}-Z_1)^{-1}\\
& \quad + \left. (Z-Z_1)(Z-\bar{Z}_1)^{-1}\left({\text d}\bar{Z}_1\right)(\bar{Z}-Z_1)^{-1}\left({\text d}Z_1\right)(\bar{Z}-Z_1)^{-1}\right\}\\
& + \left\{ - \left({\text d}Z_1\right)(Z-\bar{Z}_1)^{-1}(\bar{Z}-\bar{Z}_1)(\bar{Z}-Z_1)^{-1}\left({\text d}Z_1\right)(\bar{Z}-Z_1)^{-1} \right.\\
& \quad + (Z-Z_1)(Z-\bar{Z}_1)^{-1}\left({\text d}\bar{Z}_1\right)(Z-\bar{Z}_1)^{-1}(\bar{Z}-\bar{Z}_1)(\bar{Z}-Z_1)^{-1}\left({\text d}Z_1\right)(\bar{Z}-Z_1)^{-1} \\
& \quad +(Z-Z_1)(Z-\bar{Z}_1)^{-1}\left(-{\text d}\bar{Z}_1\right)(\bar{Z}-Z_1)^{-1}\left({\text d}Z_1\right)(\bar{Z}-Z_1)^{-1} \\
& \quad + (Z-Z_1)(Z-\bar{Z}_1)^{-1}(\bar{Z}-\bar{Z}_1)(\bar{Z}-Z_1)^{-1}\left({\text d}Z_1\right)(\bar{Z}-Z_1)^{-1}\left({\text d}Z_1\right)(\bar{Z}-Z_1)^{-1} \\
& \quad \left. + (Z-Z_1)(Z-\bar{Z}_1)^{-1}(\bar{Z}-\bar{Z}_1)(\bar{Z}-Z_1)^{-1}\left({\text d}Z_1\right)(\bar{Z}-Z_1)^{-1}\left({\text d}Z_1\right)(\bar{Z}-Z_1)^{-1} \right\}.
\endalign
$$
Thus we have
$$
2\left.d^2\frak R(Z,Z_1)\right|_{Z_1=Z}=4({\text d}Z)\left(Z-\bar Z\right)^{-1}({\text d}\bar Z)\left(\bar Z-Z\right)^{-1}
={\text d}ZY^{-1}{\text d}\bar ZY^{-1},
$$
where we have written the fixed point $Z$ as $Z=X+iY$, and we also use $Z$ to denote the variable $Z_1$ for simplification and convention of notations.
Summarizing the above, we can state the following

\proclaim{Theorem 6} The quadratic differential form ${\text
d}s^2$, which is defined by
$$
{\text d}s^2:={\text tr}\left(Y^{-1}{\text d}ZY^{-1}{\text d}\bar Z\right)={\text tr}\left({\text d}ZY^{-1}{\text d}\bar ZY^{-1}\right)
$$
at $Z$, is invariant under the action by every $M\in \hat{\Omega}_2$. And further, ${\text d}s^2$ is positive definite at any point $Z\in \hat\Bbb H_2$. So ${\text d}s^2$ defines a metric on $\hat\Bbb H_2$, it is the so-called symplectic metric on $\hat\Bbb H_2$.
\endproclaim

\demo{Proof} We only need to prove the positiveness of ${\text d}s^2$. Since the differential form ${\text d}s^2$ is invariant under the action by $M\in \hat{\Omega}_2$, and the action of $\hat{\Omega}_2$ on $\hat\Bbb H_2$ is transitive, it is sufficient to prove the positiveness of it at the special point $iI$. But there, simple computation gives
$$
{\text d}s^2={\text tr}\left(({\text d}X)^2+({\text d}Y)^2\right)=2\left({\text d}x_1^2+{\text d}x_2^2+{\text d}y_1^2+{\text d}y_2^2\right),
$$
which is plainly positive, where we have written the the corresponding variable $Z$ at the point $iI$ as $Z=X+iY=\pmatrix x_1 & x_2 \\ x_2 & x_1 \endpmatrix + i \pmatrix y_1 & y_2 \\ y_2 & y_1 \endpmatrix$. The proof of Lemma 4 is complete.
\enddemo

Next, we turn to consider the uniqueness, up to a positive constant multiple, of the symplectic metric. More precisely, we have the following

\proclaim{Lemma 4} Let $Z=X+iY=\pmatrix z_1 & z_2 \\ z_2 &
z_1\endpmatrix \in \hat\Bbb H_2$. Assume $f_{Z}\left({\text
d}z_1,\,{\text d}z_2,\,{\text d}\bar z_1,\,{\text d}\bar
z_2\right)$ is an arbitrary positive definite quadratic
differential form of ${\text d}z_1,\,{\text d}z_2,\,{\text d}\bar
z_1,\,{\text d}\bar z_2$, which is invariant under the actions of
$M\in \hat{\Omega}_2$, written as
$$
f_{Z}\left({\text d}z_1,\,{\text d}z_2,\,{\text d}\bar
z_1,\,{\text d}\bar z_2\right):= \left({\text d}z_1,\, {\text
d}z_2,\,{\text d}\bar z_1,\,{\text d}\bar
z_2\right)A\left(Z\right)\,  ^t\!\! \left({\text d}\bar
z_1,\,{\text d}\bar z_2,\,{\text d}z_1,\,{\text d}z_2\right), \tag
33
$$
where $A\left(Z\right)$ is an Hermitian matrix in $\Bbb C^{4\times 4}$ depending at most on $Z$ and $f$. Then there exists a positive constant $C$ depending only on $f$ such that
$$
f_{Z}\left({\text d}z_1,\,{\text d}z_2,\,{\text d}\bar
z_1,\,{\text d}\bar z_2\right)=C{\text tr}\left(Y^{-1}{\text
d}ZY^{-1}{\text d}\bar Z\right).\tag 34
$$
\endproclaim

\demo{Proof} First of all, by the transitivity of the action of $\hat{\Omega}_2$ on $\hat\Bbb H_2$, we can take an $M_0\in \hat{\Omega}_2$ such that $M_0<iI>=Z$. Using this $M_0$ we can make a transformation of variables $Z'=\pmatrix z'_1 & z'_2 \\ z'_2 & z'_1\endpmatrix=M_0^{-1}<Z>=M_0^{-1}\left<\pmatrix z_1 & z_2 \\ z_2 & z_1\endpmatrix\right>$. Then the special case of Lemma 5 with the point $Z=iI$ and ${\text d}Z={\text d}Z'$ gives rise to
$$
f_{iI}\left({\text d}z'_1,\,{\text d}z'_2,\,{\text d}\bar z'_1,\,{\text d}\bar z'_2\right)=C{\text tr}\left({\text d}Z'{\text d}\bar Z'\right).
$$
Using $M_0$ to act on both sides of this equality and in view of
the invariance of the both sides under this action, we deduce (34)
immediately in general. So without loss of generality we can
assume $Z=iI$ in this lemma. And thus the problem becomes to prove
(34) at the point $iI \in \hat\Bbb H_2$, i.e.,
$$
f_{iI}\left({\text d}z_1,\,{\text d}z_2,\,{\text d}\bar
z_1,\,{\text d}\bar z_2\right)=C{\text tr}\left({\text d}Z{\text
d}\bar Z\right).\tag 35
$$
Taking a step further, we consider the reflection of the problem to that of the domain $\hat E_2$. To this end, we are naturally led to use the variable transformation $Z_0=L^{-1}<Z>=(Z-iI)(Z+iI)^{-1}$. Note that, under this transformation the point $iI\in \hat\Bbb H_2$ is transformed to $0\in \hat E_2$. And we have
$$
{\text d}Z_0=({\text d}Z)(Z+iI)^{-1}-(Z-iI)(Z+iI)^{-1}({\text d}Z)(Z+iI)^{-1}.
$$
Hence, at the point $Z=iI$, or $Z_0=0$, we have
$$
{\text d}Z=2 i {\text d}Z_0. \tag 36
$$
Whence ${\text d}\bar Z=-2 i {\text d}\bar Z_0$. Consequently ${\text d}Z{\text d}\bar Z=4{\text d}Z_0{\text d}\bar Z_0$, and thus
$$
{\text Tr}\left({\text d}Z{\text d}\bar Z\right)=4{\text
Tr}\left({\text d}Z_0{\text d}\bar Z_0\right).\tag 37
$$
Further, if we write $Z_0=\pmatrix z'_1 & z'_2 \\ z'_2 &
z'_1\endpmatrix$ and recall $Z=\pmatrix z_1 & z_2 \\ z_2 &
z_1\endpmatrix$, then at the point $Z=iI$ we have by (36),
$$
\pmatrix {\text d}z_1 \\ {\text d}z_2 \endpmatrix = 2 i \pmatrix {\text d}z'_1 \\ {\text d}z'_2 \endpmatrix , \quad
\pmatrix {\text d}\bar z_1 \\ {\text d}\bar z_2 \endpmatrix =- 2 i \pmatrix {\text d}\bar z'_1 \\ {\text d}\bar z'_2 \endpmatrix .
$$
So the left hand side of (35) can be written as
$$
\align
& \quad f_{iI}\left({\text d}z_1,\,{\text d}z_2,\,{\text d}\bar z_1,\,{\text d}\bar z_2\right) \\
& = \left({\text d}z_1,\,{\text d}z_2,\,{\text d}\bar z_1,\,{\text d}\bar z_2\right) A\, ^t
\left({\text d}\bar z_1,\,{\text d}\bar z_2,\,{\text d} z_1,\,{\text d} z_2\right) \\
& =(2i)^2 \left({\text d}z'_1,\,{\text d}z'_2,\,-{\text d}\bar z'_1,\,-{\text d}\bar z'_2\right) A\, ^t
\left(-{\text d}\bar z'_1,\,-{\text d}\bar z'_2,\,{\text d} z'_1,\,{\text d} z'_2\right) \\
& =4 \left(-{\text d}z'_1,\,-{\text d}z'_2,\,{\text d}\bar z'_1,\,{\text d}\bar z'_2\right) A\, ^t
\left(-{\text d}\bar z'_1,\,-{\text d}\bar z'_2,\,{\text d} z'_1,\,{\text d} z'_2\right) \\
& =4 \left({\text d}z'_1,\,{\text d}z'_2,\,{\text d}\bar z'_1,\,{\text d}\bar z'_2\right) \pmatrix -I & 0 \\ 0 & I \endpmatrix  A  \pmatrix -I & 0 \\ 0 & I \endpmatrix \, ^t \left({\text d}\bar z'_1,\,{\text d}\bar z'_2,\,{\text d} z'_1,\,{\text d} z'_2\right) \\
& =4 \left({\text d}z'_1,\,{\text d}z'_2,\,{\text d}\bar
z'_1,\,{\text d}\bar z'_2\right)   A_0\,  ^t \left({\text d}\bar
z'_1,\,{\text d}\bar z'_2,\,{\text d} z'_1,\,{\text d} z'_2\right)
, \tag 38
\endalign
$$
where $A_0=\pmatrix -I & 0 \\ 0 & I \endpmatrix  A  \pmatrix -I &
0 \\ 0 & I \endpmatrix .$ Note that the positiveness of $A_0$ is
equivalent to that of $A$. Thus by (37) and (38), the proof of
(35) is reduced to proving
$$
\left({\text d}z'_1,\,{\text d}z'_2,\,{\text d}\bar z'_1,\,{\text d}\bar z'_2\right)  A_0\,  ^t \left({\text d}\bar z'_1,\,{\text d}\bar z'_2,\,{\text d} z'_1,\,{\text d} z'_2\right) = C {\text Tr}\left({\text d}Z_0{\text d}\bar Z_0\right),
$$
or simply, if we use again $A$ to denote $A_0$ and let $Z_0=\pmatrix z_1 & z_2 \\ z_2 & z_1 \endpmatrix$,
$$
\left({\text d}z_1,\,{\text d}z_2,\,{\text d}\bar z_1,\,{\text
d}\bar z_2\right)  A\,   ^t \left({\text d}\bar z_1,\,{\text
d}\bar z_2,\,{\text d} z_1,\,{\text d} z_2\right) = C {\text
Tr}\left({\text d}Z_0{\text d}\bar Z_0\right), \tag 39
$$
under the assumptions that the left hand side of (39), denoted by
$f_0$, is positive definite, and is invariant under the actions of
elements $M\in \hat{\Omega}_{E_2}$ with fixed point $0$. Recall
that by, Theorem 3, any action in $\hat{\Omega}_{E_2}$ with fixed
point $0$ must be of the form $Z_0 \longrightarrow \ ^tU Z_0 U$
with $U$ being unitary and satisfying $Uq=\ve qU$. Thus if we put
$ ^tU Z_0 U=Z_0'$ and let $Z_0'=\pmatrix z'_1 & z'_2 \\ z'_2 &
z'_1\endpmatrix$, then the invariance of $f_0$ is equivalent to
that the matrix $A$ satisfies
$$
\left({\text d}z'_1,\,{\text d}z'_2,\,{\text d}\bar z'_1,\,{\text
d}\bar z'_2\right)A\,  ^t\!\! \left({\text d}\bar z'_1,\,{\text
d}\bar z'_2,\, {\text d}z'_1,\,{\text d}z'_2\right)= \left({\text
d}z_1,\,{\text d}z_2,\,{\text d}\bar z_1,\,{\text d}\bar
z_2\right)A\,  ^t\!\! \left({\text d}\bar z_1,\,{\text d}\bar
z_2,\,{\text d}z_1,\,{\text d}z_2\right) \tag 40
$$
for any $Z=\pmatrix z_1 & z_2 \\ z_2 & z_1\endpmatrix$. Therefore,
our problem finally becomes to prove (39) under the assumptions of
(40) and the positiveness of it. To do this, we are naturally led
to consider at first a more explicit expression of (40). From
$Uq=\ve qU$ we can write $U=\pmatrix u_1 & u_2 \\ \ve u_2 & \ve
u_1 \endpmatrix =\pmatrix U_1 \\ \ve U_2 \endpmatrix$, where
$U_1=\pmatrix u_1, &\!\!\!\! u_2 \endpmatrix$ and $U_2=\pmatrix
u_2,&\!\!\!\! u_1 \endpmatrix$. Hence from $ Z_0'=\ ^tU Z_0 U$ we
can obtain by direct computations
$$
\align
{\text d}Z_0' & =\ ^tU {\text d}Z_0 U=\left( ^tU_1U_1+ ^tU_2U_2\right){\text d}z_1 + \ve \left( ^tU_2U_1+ ^tU_1U_2\right){\text d}z_2 \\
&= \pmatrix u_1^2 + u_2^2 & 2u_1u_2 \\ 2u_1u_2 & u_1^2 + u_2^2 \endpmatrix {\text d}z_1 + \ve \pmatrix 2u_1u_2 & u_1^2 + u_2^2 \\ u_1^2 + u_2^2 & 2u_1u_2  \endpmatrix {\text d}z_2,
\endalign
$$
whence
$$
\pmatrix {\text d}z_1' \\ {\text d}z_2' \endpmatrix = B_0 \pmatrix {\text d}z_1 \\ {\text d}z_2 \endpmatrix ,
$$
where $B_0= \pmatrix u_1^2 + u_2^2 & 2\ve u_1u_2 \\ 2u_1u_2 & \ve \left(u_1^2 + u_2^2\right) \endpmatrix .$ Thus we have
$$
\pmatrix {\text d}\bar z_1' \\ {\text d}\bar z_2' \\ {\text d}z_1'
\\ {\text d}z_2' \endpmatrix = \pmatrix \bar B_0 & 0 \\ 0 & B_0
\endpmatrix \pmatrix {\text d}\bar z_1 \\ {\text d}\bar z_2 \\
{\text d}z_1 \\ {\text d}z_2 \endpmatrix . \tag 41
$$
Denoting $B=\pmatrix \bar B_0 & 0 \\ 0 & B_0 \endpmatrix$ and
inserting (41) into (40), we see that, for any $B$ determined by
$U$ above, there always holds
$$
A=\ ^t\bar BAB.\tag 42
$$
Further, if we write $A=\left(a_{ij}\right)_{4\times 4}=\pmatrix
A_{1} & A_{2} \\ A_3 & A_4 \endpmatrix $ with $a_{ij}\in \Bbb{C}$
and $A_1,\, ...,A_4\in \Bbb{C}^{2\times 2}$, and in view of the
expanding form of the left hand side of (39), we can naturally
assume also without loss of generality that
$$
a_{33}=a_{11},\,\,a_{34}=a_{21},\,\,a_{43}=a_{12},\,\,a_{44}=a_{22},\,\,a_{23}=a_{14},\,\,a_{41}=a_{32}.
\tag 43
$$
This together with $ ^t \bar A=A$ implies
$$
 ^t \bar A_1=A_1,\,\, ^t A_2=A_2,\,\, ^t \bar A_2=A_3,\,\,^t A_1=A_4. \tag
 44
$$
Also, inserting the block expressions of $A$ and $B$ above into
(42) we obtain
$$
^t B_0A_1\bar B_0=A_1, \quad ^t B_0A_2B_0=A_2. \tag 45
$$
Again, by taking $u_1=\f{\sqrt 3}{2}i$, $u_2=\f12$, it can be
checked easily that for $\ve =\pm 1$  the matrix $U=\pmatrix u_1 &
u_2 \\ \ve u_2 & \ve u_1 \endpmatrix$ is unitary, and for this
choice of $U$, we have $u_1^2+u_2^2=-\f12$, $u_1u_2=\f{\sqrt
3}{4}i$, so $B_0=\pmatrix -\f12 & \ve \f{\sqrt 3}{2}i \\  \f{\sqrt
3}{2}i & -\f{1}{2}\ve \endpmatrix$. Inserting this $B_0$ into the
first identity of (45) we get
$$
-\f12a_{11}+\f{\sqrt 3}{2}ia_{21}=-\f12a_{11}+\ve \f{\sqrt
3}{2}ia_{12},\,\,-\f12a_{12}+\f{\sqrt 3}{2}ia_{22}=\f{\sqrt
3}{2}ia_{11}-\f12\ve a_{12}. \tag 46
$$
From the first identity of (46) we see that $a_{21}=\ve a_{12}$
holds for both of $\ve =\pm 1$, and this clearly implies
$a_{21}=0$, and whence $a_{12}=0$. And from the second equality of
(46) with $\ve =1$ we derive $a_{11}=a_{22}$. Thus we deduce
$$
A_1=a_{11}I,\tag 47
$$
and since $^t \bar A_1=A_1$ by (44), we see that $a_{11}$ is real.
Similarly, inserting the above $B_0$ into the second identity of
(45) we can get by direct computations that $-\ve a_{23}=a_{14}$
and $-\ve a_{24}=a_{13}$, and these equalities with $\ve =\pm 1$
clearly imply $a_{14}=a_{23}=0$ and $a_{13}=a_{24}=0$
respectively. Thus we have
$$
A_2=0. \tag 48
$$
The combination of (44), (47) and (48) clearly implies
$$
A=a_{11}I,
$$
and by the assumption of the positiveness of $f_0$ we obtain
$a_{11}>0.$ Substituting this $A=a_{11}I$ into the left hand side
of (39) it can be seen easily that (39) holds by taking
$C=a_{11}$. The proof of Lemma 4 is complete.
\enddemo

\head 7. The geodesics and the symplectic distances
\endhead

In this section we come to consider the geodesic and the symplectic distance of two points in $\hat\Bbb H_2$ under the metric given by the preceding section. Our first preliminary result is the following

\proclaim{Lemma 5} Under the metric ${\text d}s^2=Y^{-1}{\text
d}ZY^{-1}{\text d}\bar Z$ defined in the preceding section, the
shortest length of sufficiently smooth curves $Z=Z(t)$ joining the
points $Z(0)=iI$ and $Z(t_0)=i\L=i\pmatrix \lambda_1 & \lambda_2\\
\lambda_2 & \lambda_1 \endpmatrix$ with $\lambda_1\geq\lambda_2+1$
and $\lambda_2\geq 0$ is
$$
\rho=\left(\log^2(\lambda_1+\lambda_2)+\log^2(\lambda_1-\lambda_2)\right)^{1/2},
$$
where $t_0$ is a positive constant and $0\le t\le t_0$ is a parameter. And a curve of shortest length can be taken as
$$
Z=Z(t)=i\pmatrix\frac{(\lambda_1+\lambda_2)^{t/t_0}+(\lambda_1-\lambda_2)^{t/t_0}}{2} & \frac{(\lambda_1+\lambda_2)^{t/t_0}-(\lambda_1-\lambda_2)^{t/t_0}}{2}\\
\frac{(\lambda_1+\lambda_2)^{t/t_0}-(\lambda_1-\lambda_2)^{t/t_0}}{2}
&
\frac{(\lambda_1+\lambda_2)^{t/t_0}+(\lambda_1-\lambda_2)^{t/t_0}}{2}\endpmatrix
. \tag 49
$$
\endproclaim

\demo{Proof} For any curve $Z=Z(t)$ under consideration we write
$Z(t)=X(t)+iY(t)$. Then we have
$$
X(0)=0,\,\,Y(0)=I,\,\,X(t_0)=0,\,\,Y(t_0)=\L . \tag 50
$$
The length $S$ of this curve then is
$$
S=\int_{0}^{t_0}\left( tr\!\! \left( Y(t)^{-1}\dot Z(t)Y(t)^{-1}\dot{ \overline {Z(t)}}\right) \right)^{1/2}{\text d}t ,
$$
where $\dot Z(t)$ is defined as $\dot Z(t)={\text d}Z(t)/{\text d}t$. Now direct computation yields
$$
tr\!\! \left( Y(t)^{-1}\dot Z(t)Y(t)^{-1}\dot{ \overline {Z(t)}}\right)
=tr\!\! \left( Y(t)^{-1}\dot X(t)Y(t)^{-1}\dot X(t)\right)+tr\!\! \left( Y(t)^{-1}\dot Y(t)Y(t)^{-1}\dot Y(t)\right),
$$
and
$$
tr\!\! \left( Y(t)^{-1}\dot X(t)Y(t)^{-1}\dot X(t)\right)\ge 0.
$$
Thus we have
$$
S\ge \int_{0}^{t_0}\left( tr\!\! \left( Y(t)^{-1}\dot
Y(t)Y(t)^{-1}\dot Y(t)\right)\right)^{1/2}{\text d}t , \tag 51
$$
and the equality holds if $\dot X(t)=0$. Further, if we write $Y(t)=\pmatrix y_1(t) & y_2(t)\\ y_2(t) & y_1(t) \endpmatrix$ with $y_1(t)>\left|y_2(t)\right|$ then by direct computations we have
$$
\align
S_1&:=\int_{0}^{t_0}\left( tr\!\! \left( Y(t)^{-1}\dot Y(t)Y(t)^{-1}\dot Y(t)\right)\right)^{1/2}{\text d}t\\
&=\int_0^{t_0}\left((\dot y_1(t)+\dot
y_2(t))^2(y_1(t)+y_2(t))^{-2}+(\dot y_1(t)-\dot
y_2(t))^2(y_1(t)-y_2(t))^{-2}\right)^{1/2}{\text d}t.\tag 52
\endalign
$$
Put
$$
c_1=\rho^{-1}\log(\lambda_1+\lambda_2),\ c_2=\rho^{-1}\log(\lambda_1-\lambda_2).
$$
Then we have clearly $c_1,\,c_2\geq 0,$ $c_1^2+c_2^2=1$, and
$$
\left(\frac{\dot{y_1}+\dot{y_2}}{y_1+y_2}\right)^2+\left(\frac{\dot{y_1}-\dot{y_2}}{y_1-y_2}\right)^2
=\left(c_1\frac{\dot{y_1}+\dot{y_2}}{y_1+y_2}+c_2\frac{\dot{y_1}-\dot{y_2}}{y_1-y_2}\right)^2
+\left(c_1\frac{\dot{y_1}-\dot{y_2}}{y_1-y_2}-c_2\frac{\dot{y_1}+\dot{y_2}}{y_1+y_2}\right)^2 .
$$
So
$$
S_1\geq\int_0^{t_0}\left|c_1\frac{\dot{y_1}+\dot{y_2}}{y_1+y_2}+c_2\frac{\dot{y_1}-\dot{y_2}}{y_1-y_2}\right|{\text d}t ,
$$
and
$$
S_1=\int_0^{t_0}\left(c_1\frac{\dot{y_1}+\dot{y_2}}{y_1+y_2}+c_2\frac{\dot{y_1}-\dot{y_2}}{y_1-y_2}\right){\text
d}t \tag 53
$$
if and only if
$$
\cases
&\!\!\!\!\! c_1\frac{\dot{y_1}+\dot{y_2}}{y_1+y_2}+c_2\frac{\dot{y_1}-\dot{y_2}}{y_1-y_2}\geq0,\\
&\!\!\!\!\! {c_1\frac{\dot{y_1}-\dot{y_2}}{y_1-y_2}-c_2\frac{\dot{y_1}+\dot{y_2}}{y_1+y_2}=0}.
\endcases
\tag 54
$$
Thus by (53) and (50) we have under (54),
$$
\align
S_1=&\int_0^{t_0}\left(c_1\frac{\dot{y_1}+\dot{y_2}}{y_1+y_2}+c_2\frac{\dot{y_1}-\dot{y_2}}{y_1-y_2}\right)dt\\
=&\left(c_1\log(y_1+y_2)+c_2\log(y_1-y_2)\right)\mid_0^1\\
=&c_1\log(\lambda_1+\lambda_2)+c_2\log(\lambda_1-\lambda_2)=\rho .
\tag 55
\endalign
$$
This proves that the shortest length of sufficiently smooth curves from the point $Z(0)=iI$ to the point $Z(t_0)=i\Lambda$ to be equal at least to
$\rho$.
Again, if we take a curve by choosing
$$
Z(t)=i\pmatrix\frac{(\lambda_1+\lambda_2)^{t/t_0}+(\lambda_1-\lambda_2)^{t/t_0}}{2} & \frac{(\lambda_1+\lambda_2)^{t/t_0}-(\lambda_1-\lambda_2)^{t/t_0}}{2}\\
\frac{(\lambda_1+\lambda_2)^{t/t_0}-(\lambda_1-\lambda_2)^{t/t_0}}{2} & \frac{(\lambda_1+\lambda_2)^{t/t_0}+(\lambda_1-\lambda_2)^{t/t_0}}{2}\endpmatrix ,
$$
then clearly $Z(0)=iI,\ Z(t_0)=i\Lambda,\ X(t)=0$ and it can be
checked by direct computations easily that all the conditions of
(54) are satisfied. Therefore the length of this curve is indeed
equal exactly to
$\rho=\sqrt{\log^2(\lambda_1+\lambda_2)+\log^2(\lambda_1-\lambda_2)}$,
and the proof of this lemma is complete.
\enddemo

The next lemma deals with the differential equation of the
geodesics and the uniqueness of the curve with shortest length in
the previous lemma.

\proclaim{Lemma 6} Under the notations of Lemma 5, {\text (49)} is
the unique curve of shortest length $\rho$, and it is called the
geodesic through the points $iI$ and $i\L$. Furthermore, the
parameter $t$ is nothing but the the length of the arcs in the
geodesic starting from the point $Z(0)=iI$. The differential
equation of the geodesic $Z=Z(s)$ is
$$
\ddot{Z}(s)=-i\dot ZY^{-1}\dot Z.\tag 56
$$
\endproclaim

\demo{Proof} Assume that $Z=Z(s)$ is an arbitrary piecewise
sufficiently smooth curve connecting the points $Z(0)=iI$ and
$Z\left(s_0\right)=i\L$, where the parameter $s$ denotes of the
length of the arc. Then by the definition ${\text
d}s^2=tr\left(Y^{-1}{\text d}ZY^{-1}{\text d}\bar Z\right)$ of our
metric we have
$$
tr\left(Y(s)^{-1}\dot{Z}(s)Y(s)^{-1}\dot{\bar{Z}}(s)\right)=tr\left(Y(s)^{-1}\f{{\text d}Z(s)}{ds}Y(s)^{-1}\f{{\text d}\bar{Z}(s)}{ds}\ri)=\f{ds^2}{ds^2}=1,
$$
where $\dot{Z}(s)=\f{{\text d}Z(s)}{ds}$, $\dot{\bar{Z}}(s)=\f{{\text d}\bar{Z}(s)}{ds}$. Thus the length $s_0$ of this curve can be written as
$\int_0^{s_0}tr\left(Y(s)^{-1}\dot{Z}(s)Y(s)^{-1}\dot{\bar{Z}}(s)\right)^{1/2}{\text d}s$, or abbreviated as
$$
\int_0^{s_0}tr\left(Y^{-1}\dot{Z}Y^{-1}\dot{\bar{Z}}\right)^{1/2}{\text
d}s .\tag 57
$$
Now let $Z(s)$ vary with the same endpoints $Z(0)=iI$ and
$Z\left(s_0\right)=i\L$ as the given curve of the lemma, then (57)
can be viewed as a functional of the varied curves $Z=Z(s)$.
Moreover, if we use the parameter $s_1$ to denote the length of
the arcs of the varied curves $Z=Z(s)$, and use $S_1$ to denote
the length of this curve, then by changing the integral variable
$s$ to $s_1$ we see that (57), with $Z$ to be the present curve
$Z=Z(s)$, is equal to
$$
\int_0^{s_0}tr\left(Y^{-1}\f{{\text d} Z}{{\text d} s}Y^{-1}\f{{\text d} \bar Z}{{\text d} s}\right)^{1/2}{\text d}s
= \int_0^{S_1}tr\left(Y^{-1}\f{{\text d} Z}{{\text d} s_1}Y^{-1}\f{{\text d} \bar Z}{{\text d} s_1}\right)^{1/2}{\text d}s_1
= \int_0^{S_1}{\text d}s_1=S_1 .
$$
This shows that (57) is always equal to the length of the curves
with endpoints $Z(0)=iI$ and $Z\left(s_0\right)=i\L$ whether or
not the parameter $s$ denotes the length of arc in the curve. And
hence, as a functional of the curve $Z=Z(s)$ connecting the points
$Z(0)=iI$ and $Z\left(s_0\right)=i\L$, (57) is indeed the length
of the curves. Therefore, for any curve $Z=Z(s)$ with endpoints
$Z(0)=iI$ and $Z\left(s_0\right)=i\L$ of shortest length $\rho$,
and with the parameter $s$ to be the length of arc of this curve,
we have
$$
\delta
\int_0^{s_0}tr\left(Y^{-1}\dot{Z}Y^{-1}\dot{\bar{Z}}\right)^{1/2}{\text
d}s =0, \tag 58
$$
where $\d J(Z(s))$ denotes the first variation of a functional
$J(Z(s))$ at $Z(s)$. Further, for any function $f(s)$, let $\d
f(s)$ denote the first variation of $f(s)$, and let $\d Z(s)$
denote the matrix with the corresponding entries replaced by the
first variations of them. Then by the simple properties of the
first variation we get from (58),
$$
\int_0^{s_0}\delta
\left(tr\left(Y^{-1}\dot{Z}Y^{-1}\dot{\bar{Z}}\right)\right){\text
d}s =0. \tag 59
$$
Set $B=\dot{Z}Y^{-1}\dot{\bar{Z}},\ W=Y^{-1}\dot{Z}Y^{-1}$. Then we have by direct computations,
$$
\align
&\delta \left(tr\left(Y^{-1}\dot{Z}Y^{-1}\dot{\bar{Z}}\right)\right) =  tr\left(\d\left(Y^{-1}\dot{Z}Y^{-1}\dot{\bar{Z}}\right)\right)\\
=&tr\left(\left(\d Y^{-1}\right)\dot{Z}Y^{-1}\dot{\bar{Z}}+Y^{-1}\d(\dot{Z})Y^{-1}\dot{\bar{Z}}+Y^{-1}\dot{Z}(\d Y^{-1})\dot{\bar{Z}}
+Y^{-1}\dot{Z}Y^{-1}\d(\dot{\bar{Z}})\right)\\
=&2{\text {Re}} tr\left(B(\d Y^{-1})+(W\d
\bar{Z})^{.}-\dot{W}\d\bar{Z}\right).\tag 60
\endalign
$$
Again from $YY^{-1}=I$ we can derive easily $(\d Y)Y^{-1}+Y\d Y^{-1}=0$ and $\left(Y^{-1}\right)^{.}Y+Y^{-1}\dot Y=0$, hence we have
$$
\d Y^{-1}=-Y^{-1}(\d Y)Y^{-1}=\f{i}{2}Y^{-1}\left(\d Z-\d \bar
Z\right)Y^{-1},\tag 61
$$
and
$$
\left(Y^{-1}\right)^{.}=-Y^{-1}\dot YY^{-1}=\f{i}{2}Y^{-1}\left(\dot Z-\dot{\bar Z}\right)Y^{-1}.
$$
And by the last equality we have
$$
\align
\dot W & =\left(Y^{-1}\dot ZY^{-1}\right)^{.}=Y^{-1}\ddot ZY^{-1}+\left(Y^{-1}\right)^{.}\dot ZY^{-1}+Y^{-1}\dot Z \left(Y^{-1}\right)^{.}\\
& =Y^{-1}\ddot ZY^{-1}+iY^{-1}\dot ZY^{-1}\dot
ZY^{-1}-\f{i}2Y^{-1}\left(B+\bar B\right)Y^{-1}.\tag 62
\endalign
$$
Inserting (61) and (62) into the last expression of (60), by
direct computations, we can write (60) further as
$$
=2{\text {Re}} tr\left((W\d \bar{Z})^{.}-Y^{-1}\ddot ZY^{-1}\d {\bar Z}-iY^{-1}\dot ZY^{-1}\dot ZY^{-1}\d {\bar Z}\right).
$$
Inserting this into (59) and in view of
$$
\int_0^{s_0}{\text {Re tr}}\left((W\d \bar{Z})^{.}\right) {\text
d}s=0
$$
since $\d Z=0$ at the end points with $s=0$ and $s=s_0$, we get
$$
{\text {Re}}\int_0^{s_0}{\text tr}\left(Y^{-1}\left(\ddot{Z}+i\dot{Z}Y^{-1}\dot{Z}\right)Y^{-1}\d \bar{Z}\right){\text d}s=0.
$$
This clearly implies
$$
\ddot{Z}+i\dot{Z}Y^{-1}\dot{Z}=0,\tag 63
$$
which is the differential equation satisfied by the geodesics. Now
we are in the position to consider the solution of this
differential equation under the boundary conditions $Z(0)=iI$ and
$Z(s_0)=i\L$. In view of $Z(s)=\pmatrix z_1(s) & z_2(s) \\ z_2(s)
& z_1(s)\endpmatrix$, we can use $p^{-1}=\f{1}{\sqrt 2}\pmatrix
1&1\\ -1&1\endpmatrix$ and $p=\f{1}{\sqrt 2}\pmatrix 1&-1\\
1&1\endpmatrix$ to multiply the left hand side of (63) from the
left and the right respectively and obtain
$$
\pmatrix \ddot{{z}}_1+\ddot{{z}}_2 & 0\\ 0 & \ddot{{z}}_1-\ddot{{z}}_2 \endpmatrix
=-i\pmatrix \dot{z_1}+\dot{z_2} & 0\\ 0 & \dot{z_1}-\dot{z_2} \endpmatrix
\pmatrix (y_1+y_2)^{-1} & 0\\ 0 & (y_1-y_2)^{-1} \endpmatrix
\pmatrix \dot{z_1}+\dot{z_2} & 0\\ 0 & \dot{z_1}-\dot{z_2} \endpmatrix .
$$
Here we have used the abbreviations $z_1=z_1(s)$, $z_2=z_2(s)$, and the symbols $z_1=x_1+iy_1=x_1(s)+iy_1(s)$, $z_2=x_2+iy_2=x_2(s)+iy_2(s)$, with real $x_1,\ y_1,\ x_2,\ y_2$. Hence
$$
\ddot{{z}}_1+\ddot{{z}}_2=-i\left(\dot{z_1}+\dot{z_2}\right)\left(y_1+y_2\right)^{-1},\quad
\ddot{{z}}_1-\ddot{{z}}_2=-i\left(\dot{z_1}-\dot{z_2}\right)\left(y_1-y_2\right)^{-1},
$$
or equivalently, by considering the real and the imaginary parts respectively,
$$
\ddot{{x}}_1+\ddot{{x}}_2=\f{2\left(\dot x_1+\dot
x_2\right)\left(\dot y_1+\dot
y_2\right)}{\left(y_1+y_2\right)},\,\,\,\,
\ddot{{y}}_1+\ddot{{y}}_2=\f{\left(\dot{y_1}+\dot{y_2}\right)^2-\left(\dot{x_1}+\dot{x_2}\right)^2}{\left(y_1+y_2\right)},
\tag 64
$$
and
$$
\ddot{{x}}_1-\ddot{{x}}_2=\f{2\left(\dot{x_1}-\dot{x_2}\right)\left(\dot{y_1}-\dot{y_2}\right)}{\left(y_1-y_2\right)},\,\,\,\,
\ddot{{y}}_1-\ddot{{y}}_2=\f{\left(\dot{y_1}-\dot{y_2}\right)^2-\left(\dot{x_1}-\dot{x_2}\right)^2}{\left(y_1-y_2\right)}.
\tag 65
$$
By the first equality of (64) one can derive easily that
$$
\dot x_1+\dot x_2 = c_1\left({y_1}+{y_2}\right)^2,
$$
where $c_1$ is a constant. This together with the boundary conditions $x_1(0)=x_2(0)=0$ implies, for any $s$ with $0\le s\le s_0$,
$$
x_1(s)+x_2(s)=c_1\int_0^s \left({y_1}(s)+{y_2}(s)\right)^2{\text d}s.
$$
This in combination with the boundary conditions $x_1(s_0)=x_2(s_0)=0$ implies $c_1=0$, whence we get
$$
x_1(s)+x_2(s)=0.
$$
Similarly, by the first equality of (65) one can derive
$$
x_1(s)-x_2(s)=0.
$$
Thus we have
$$
x_1(s)=x_2(s)=0.\tag 66
$$
Again, substituting (66) into the second equality of (64) we see
that
$$
\ddot{{y}}_1+\ddot{{y}}_2={\left(\dot y_1+\dot y_2\right)^2}{\left(y_1+y_2\right)}^{-1}.
$$
This clearly implies $\dot y_1+\dot y_2 = c_2 \left(y_1+y_2\right),$ and so
$$
y_1+y_2=e^{c_2s+c_3},
$$
where $c_2$ and $c_3$ again are constants. This together with the boundary conditions $y_1(0)=1,\,\,y_2(0)=0,\,\,y_1(s_0)=\l_1$ and $y_2(s_0)=\l_2$ gives rise to
$$
e^{c_3}=1,\,\,\,e^{c_2s_0+c_3}=\l_1+\l_2,
$$
hence $c_3=0$, $c_2=s_0^{-1}\log \left(\l_1+\l_2\right)$. Thus
$$
y_1(s)+y_2(s)=\left(\l_1+\l_2\right)^{s/s_0}.\tag 67
$$
Similarly from the second equality of (65) we can obtain
$$
y_1(s)-y_2(s)=\left(\l_1-\l_2\right)^{s/s_0}.\tag 68
$$
Finally, by (66), (67) and (68) we see that the solution of (63)
with the boundary conditions $Z(0)=iI$ and $Z(s_0)=i\L$ is
$$
Z=Z(s)=i\pmatrix\frac{(\lambda_1+\lambda_2)^{s/s_0}+(\lambda_1-\lambda_2)^{s/s_0}}{2} & \frac{(\lambda_1+\lambda_2)^{s/s_0}-(\lambda_1-\lambda_2)^{s/s_0}}{2}\\
\frac{(\lambda_1+\lambda_2)^{s/s_0}-(\lambda_1-\lambda_2)^{s/s_0}}{2}
&
\frac{(\lambda_1+\lambda_2)^{s/s_0}+(\lambda_1-\lambda_2)^{s/s_0}}{2}\endpmatrix
.
$$
The proof of this lemma is complete.
\enddemo

Next, we take a step further to consider the distance and the geodesic equation connecting two arbitrary points in $\hat\Bbb H_2$. To this end, we will first need a further understanding of the action of $\hat\O_2$ on $\hat\Bbb H_2$. This is the following

\proclaim{Lemma 7} Let $M=\pmatrix A & B \\ C & D \endpmatrix$ be
an element of $\hat\O_2$ with $A=\pmatrix a_1 & a_2 \\ \ve a_2
&\ve a_1 \endpmatrix$, $B=\pmatrix b_1 & b_2 \\ \ve b_2 &\ve b_1
\endpmatrix$, $C=\pmatrix c_1 & c_2 \\ \ve c_2 &\ve c_1
\endpmatrix$ and $D=\pmatrix d_1 & d_2 \\ \ve d_2 &\ve d_1
\endpmatrix$. Put
$$
M_1=\pmatrix a_1+a_2 & b_1+b_2 \\ c_1+c_2 & d_1+d_2 \endpmatrix ,\quad M_2=\pmatrix a_1-a_2 & b_1-b_2 \\ c_1-c_2 & d_1-d_2 \endpmatrix .
$$
Then we have $M_1,\,M_2\in SL_2 (\Bbb R)$, and the action of $M$ on an element $Z=\pmatrix z_1 & z_2 \\ z_2 & z_1 \endpmatrix $ of $\hat\Bbb H_2$ can be expressed as
$$
M<Z>=\pmatrix \f{M_1<z_1+z_2>+M_2<z_1-z_2>}{2} & \f{\ve
\left(M_1<z_1+z_2>-M_2<z_1-z_2>\right)}{2}
 \\ \f{\ve \left(M_1<z_1+z_2>-M_2<z_1-z_2>\right)}{2} & \f{M_1<z_1+z_2>+M_2<z_1-z_2>}{2} \endpmatrix ,
$$
where the action $M<z>$ of an element $M$ in $SL_2 (\Bbb R)$ on a point $z$ in the upper half plane is defined by the linear fractional transformation as usual.
\endproclaim

\demo{Proof} Recall that $p=\f1{\sqrt 2}\pmatrix 1 & -1 \\ 1 & 1 \endpmatrix$, $ ^t p=p^{-1}=\f1{\sqrt 2}\pmatrix 1 & 1 \\ -1 & 1 \endpmatrix$. So direct computation yields
$$
\align
^t \left( ^t pDp\right)\, ^t pAp & = \pmatrix \f{\ve +1}{2}\left(d_1+d_2\right) & \f{\ve -1}{2}\left(d_1+d_2\right) \\ \f{\ve -1}{2}\left(d_1-d_2\right) &
\f{\ve +1}{2}\left(d_1-d_2\right) \endpmatrix
\pmatrix \f{\ve +1}{2}\left(a_1+a_2\right) & \f{\ve -1}{2}\left(a_1-a_2\right) \\ \f{\ve -1}{2}\left(a_1+a_2\right) &
\f{\ve +1}{2}\left(a_1-a_2\right) \endpmatrix \\
& = \pmatrix \left(a_1+a_2\right)\left(d_1+d_2\right) & 0 \\ 0 & \left(a_1-a_2\right)\left(d_1-d_2\right) \endpmatrix ,
\endalign
$$
and similarly
$$
^t \left( ^t pBp\right)\, ^t pCp = \pmatrix \left(b_1+b_2\right)\left(c_1+c_2\right) & 0 \\ 0 & \left(b_1-b_2\right)\left(c_1-c_2\right) \endpmatrix .
$$
Thus by $ ^t DA-\, ^tBC=I$ we get
$$
\align
I & =\ ^tp ^tDp ^tpAp-\ ^tp ^tBp ^tpCp \\
& = \pmatrix \left(a_1+a_2\right)\left(d_1+d_2\right) - \left(b_1+b_2\right)\left(c_1+c_2\right) & 0 \\ 0 & \left(a_1-a_2\right)\left(d_1-d_2\right)-\left(b_1-b_2\right)\left(c_1-c_2\right) \endpmatrix .
\endalign
$$
This proves
$$
\left(a_1+a_2\right)\left(d_1+d_2\right) - \left(b_1+b_2\right)\left(c_1+c_2\right)=1
$$
and
$$
\left(a_1-a_2\right)\left(d_1-d_2\right)-\left(b_1-b_2\right)\left(c_1-c_2\right)=1 .
$$
That is to say $M_1$ and $M_2$ are in $SL_2(\Bbb R)$. Again, from the definition $M<Z>=(AZ+B)(CZ+D)^{-1}$ we can obtain by direct computations,
$$
\align
& \quad ^tpM<Z>p = \left( ^tpAp ^tpZp+ ^tpBp\right)\left( ^tpCp ^tpZp+ ^tpDp\right)^{-1} \\
& = \left( \pmatrix \f{\ve +1}{2}\left(a_1+a_2\right) & \f{\ve -1}{2}\left(a_1-a_2\right) \\ \f{\ve -1}{2}\left(a_1+a_2\right) &
\f{\ve +1}{2}\left(a_1-a_2\right) \endpmatrix
\pmatrix z_1+z_2 & 0 \\ 0 & z_1-z_2 \endpmatrix
+ \pmatrix \f{\ve +1}{2}\left(b_1+b_2\right) & \f{\ve -1}{2}\left(b_1-b_2\right) \\ \f{\ve -1}{2}\left(b_1+b_2\right) &
\f{\ve +1}{2}\left(b_1-b_2\right) \endpmatrix \right) \\
& \times \left( \pmatrix \f{\ve +1}{2}\left(c_1+c_2\right) &
\f{\ve -1}{2}\left(c_1-c_2\right) \\ \f{\ve
-1}{2}\left(c_1+c_2\right) & \f{\ve +1}{2}\left(c_1-c_2\right)
\endpmatrix \pmatrix z_1+z_2 & 0
\\ 0 & z_1-z_2 \endpmatrix + \pmatrix \f{\ve
+1}{2}\left(d_1+d_2\right) & \f{\ve -1}{2}\left(d_1-d_2\right) \\
\f{\ve -1}{2}\left(d_1+d_2\right) &
\f{\ve +1}{2}\left(d_1-d_2\right) \endpmatrix \right)^{-1} \\
& = \pmatrix \f{\ve +1}{2}\left(\left(a_1+a_2\right)\left(z_1+z_2\right)+b_1+b_2\right) & \f{\ve -1}{2}\left(\left(a_1-a_2\right)\left(z_1-z_2\right)+b_1-b_2\right) \\ \f{\ve -1}{2}\left(\left(a_1+a_2\right)\left(z_1+z_2\right)+b_1+b_2\right) &
\f{\ve +1}{2}\left(\left(a_1-a_2\right)\left(z_1-z_2\right)+b_1-b_2\right) \endpmatrix \\
& \quad \times  \pmatrix \f{\ve +1}{2}\left(\left(c_1+c_2\right)\left(z_1+z_2\right)+d_1+d_2\right) & \f{\ve -1}{2}\left(\left(c_1-c_2\right)\left(z_1-z_2\right)+d_1-d_2\right) \\ \f{\ve -1}{2}\left(\left(c_1+c_2\right)\left(z_1+z_2\right)+d_1+d_2\right) &
\f{\ve +1}{2}\left(\left(c_1-c_2\right)\left(z_1-z_2\right)+d_1-d_2\right) \endpmatrix ^{-1} .
\endalign
$$
When $\ve =1$, this can be computed further as
$$
= \pmatrix M_1<z_1+z_2> & 0 \\ 0 & M_2<z_1-z_2> \endpmatrix ,
$$
and when $\ve =-1$, it can be computed as
$$
= \pmatrix M_2<z_1-z_2> & 0 \\ 0 & M_1<z_1+z_2> \endpmatrix .
$$
So in summary we always have
$$
M<Z>=\pmatrix \f{M_1<z_1+z_2>+M_2<z_1-z_2>}{2} & \f{\ve \left(M_1<z_1+z_2>-M_2<z_1-z_2>\right)}{2} \\ \f{\ve \left(M_1<z_1+z_2>-M_2<z_1-z_2>\right)}{2} & \f{M_1<z_1+z_2>+M_2<z_1-z_2>}{2} \endpmatrix ,
$$
as desired. The proof of this lemma is complete.
\enddemo

The following elementary result will be useful for our further arguments.

\proclaim{Lemma 8} Let $\l_1,\,\l_2,\,\mu_1,\,\mu_2,\,\theta_1,\,
\theta_2$ be real numbers satisfying $\l_1\ne 0,\,\l_2\ne 0$.
Suppose $\l$ with $\l\ne 0$ is a real number satisfying the
equality
$$
\pmatrix \cos \theta_1 & \sin \theta_1 \\ -\sin \theta_1 & \cos
\theta_1 \endpmatrix \pmatrix \l_1 & \mu_1 \\ 0 & \l_1^{-1}
\endpmatrix = \pmatrix \l & 0 \\ 0 & \l^{-1} \endpmatrix \pmatrix
\cos \theta_2 & \sin \theta_2 \\ -\sin \theta_2 & \cos \theta_2
\endpmatrix \pmatrix \l_2 & \mu_2 \\ 0 & \l_2^{-1} \endpmatrix .
\tag 69
$$
Then we have
$$
\l^2+\l^{-2}=\left(\l_2/\l_1\right)^2+\left(\l_1/\l_2\right)^2+\left(\l_1\mu_2-\l_2\mu_1\right)^2.
$$
\endproclaim

\demo{Proof} From (69) we have by direct computations
$$
\align
\pmatrix \cos \theta_1 & \sin \theta_1 \\ -\sin \theta_1 & \cos \theta_1 \endpmatrix
& = \pmatrix \l & 0 \\ 0 & \l^{-1} \endpmatrix
\pmatrix \cos \theta_2 & \sin \theta_2 \\ -\sin \theta_2 & \cos \theta_2 \endpmatrix
\pmatrix \l_2 & \mu_2 \\ 0 & \l_2^{-1} \endpmatrix
\pmatrix \l_1^{-1} & -\mu_1 \\ 0 & \l_1 \endpmatrix \\
& = \pmatrix \l & 0 \\ 0 & \l^{-1} \endpmatrix
\pmatrix \cos \theta_2 & \sin \theta_2 \\ -\sin \theta_2 & \cos \theta_2 \endpmatrix
\pmatrix \l_2\l_1^{-1} & \l_1\mu_2-\l_2\mu_1 \\ 0 & \l_1\l_2^{-1} \endpmatrix \\
& = \pmatrix \l\l_2\l_1^{-1}\cos \theta_2 & \l \left(\l_1\mu_2-\l_2\mu_1\right) \cos \theta_2 + \l\l_1\l_2^{-1}\sin \theta_2 \\
-\l^{-1}\l_2\l_1^{-1}\sin \theta_2 & \l^{-1}\l_1\l_2^{-1}\cos
\theta_2 - \l^{-1}\left(\l_1\mu_2-\l_2\mu_1\right)\sin \theta_2
\endpmatrix .\\ \tag 70
\endalign
$$
By the first column of (70) we have $\pmatrix  \cos\theta_1 \\
-\sin\theta_1 \endpmatrix = \pmatrix \l\l_2\l_1^{-1}\cos \theta_2
\\ -\l^{-1}\l_2\l_1^{-1}\sin \theta_2 \endpmatrix$, thus
$$
\pmatrix  \cos\theta_1 \\ \sin\theta_1 \endpmatrix = \pmatrix
\l\l_2\l_1^{-1} & 0 \\ 0 & \l^{-1}\l_2\l_1^{-1} \endpmatrix
\pmatrix  \cos\theta_2 \\ \sin\theta_2 \endpmatrix . \tag 71
$$
And by the second column of (70) we have $\pmatrix \sin\th_1 \\
\cos\th_1 \endpmatrix = \pmatrix \l
\left(\l_1\mu_2-\l_2\mu_1\right) \cos \theta_2 +
\l\l_1\l_2^{-1}\sin \theta_2 \\ \l^{-1}\l_1\l_2^{-1}\cos \theta_2
- \l^{-1}\left(\l_1\mu_2-\l_2\mu_1\right)\sin \theta_2
\endpmatrix$, so
$$
\pmatrix  \cos\theta_1 \\ \sin\theta_1 \endpmatrix = \pmatrix
\l^{-1}\l_1\l_2^{-1} & - \l^{-1}\left(\l_1\mu_2-\l_2\mu_1\right)
\\ \l \left(\l_1\mu_2-\l_2\mu_1\right) & \l\l_1\l_2^{-1}
\endpmatrix \pmatrix  \cos\theta_2 \\ \sin\theta_2 \endpmatrix .
\tag 72
$$
Consider the difference of (71) and (72), getting
$$
 \pmatrix \l\l_2\l_1^{-1} - \l^{-1}\l_1\l_2^{-1} & \l^{-1}\left(\l_1\mu_2-\l_2\mu_1\right) \\ -\l \left(\l_1\mu_2-\l_2\mu_1\right) & \l^{-1}\l_2\l_1^{-1} - \l\l_1\l_2^{-1} \endpmatrix
\pmatrix  \cos\theta_2 \\ \sin\theta_2 \endpmatrix =0.
$$
Now in view of $\pmatrix  \cos\theta_2 \\ \sin\theta_2 \endpmatrix \ne 0$ we deduce
$$
\det \pmatrix \l\l_2\l_1^{-1} - \l^{-1}\l_1\l_2^{-1} & \l^{-1}\left(\l_1\mu_2-\l_2\mu_1\right) \\ -\l \left(\l_1\mu_2-\l_2\mu_1\right) & \l^{-1}\l_2\l_1^{-1} - \l\l_1\l_2^{-1} \endpmatrix =0,
$$
which implies the desired result. The proof of this lemma is complete.
\enddemo

\proclaim{Lemma 9} For any $z=x+iy$ in the upper half plane and
any real $\mu >0$, let $M\in SL_2(\Bbb R)$ be such that $M<z>=\mu
i.$ Then $M$ must be of the form
$$
M=\pmatrix \mu^{1/2} & 0 \\ 0
& \mu^{-1/2} \endpmatrix \pmatrix \cos\th & \sin\th \\ -\sin\th &
\cos\th \endpmatrix \pmatrix y^{-1/2} & -xy^{-1/2} \\ 0 & y^{1/2}
\endpmatrix ,
$$
where $\th$ is a real number.
\endproclaim

\demo{Proof} Let $M=\pmatrix a & b \\ c & d \endpmatrix$, then $M^{-1}=\pmatrix d & -b \\ -c & a \endpmatrix$. Thus from $M<z>=\mu i$ we get
$$
z=M^{-1}<\mu i>=\f{d\mu i-b}{ -c\mu i+a}=\f{(a+c\mu i)(d\mu i-b)}{a^2+(c\mu )^2}=-\f{ab+cd\mu^2}{a^2+(c\mu )^2}+i\f{\mu}{a^2+(c\mu )^2}.
$$
Hence we have
$$
x=-\f{ab+cd\mu^2}{a^2+(c\mu )^2},\,\,y=\f{\mu}{a^2+(c\mu )^2}.\tag
73
$$
From the latter of (73) we see that $|-c\mu
i+a|=\left(\f{\mu}{y}\right)^{1/2}$, thus we have
$$
-c\mu i+a=|-c\mu i+a|e^{i\th}=\left(\f{\mu}{y}\right)^{1/2}\cos\th +i\left(\f{\mu}{y}\right)^{1/2}\sin\th ,
$$
where $\th $ is a real number satisfying $0\le \th <2\pi$. This clearly implies
$$
a=\left(\mu y^{-1}\right)^{1/2}\cos\th , \,\,\,c=-(\mu
y)^{-1/2}\sin\th .\tag 74
$$
Again from $x+iy=\f{d\mu i-b}{ -c\mu i+a}$ we get
$$
d\mu i-b=(x+iy)(a-c\mu i)=ax+c\mu i+i(ay-c\mu x),
$$
and hence
$$
\align
& b=-ax-c\mu y=-\mu^{1/2}xy^{-1/2}\cos\th +\mu^{1/2}y^{1/2}\sin\th ,\\
& d=\mu^{-1}ay-cx=\mu^{-1/2}y^{1/2}\cos\th +\mu^{-1/2}xy^{-1/2}\sin\th .
\endalign
$$
This together with (74) completes the proof of this lemma.
\enddemo

\proclaim{Lemma 10} Let $z_1=x_1+iy_1$ and $z_2=x_2+iy_2$ be two
arbitrary points in the upper half plane, then there exists an
$M\in SL_2(\Bbb R)$ such that
$$
M<z_1>=i,\,\,\,M<z_2>=i\l ,
$$
with $\l \ge 1$ determined uniquely by $z_1$ and $z_2$. And in
this case, $\l$ satisfies
$$
\l + \l^{-1}=y_1y_2^{-1}+y_2y_1^{-1}+\left(x_1-x_2\right)^2y_1^{-1}y_2^{-1}.
$$
\endproclaim

\demo{Proof} Since the action of $SL_2(\Bbb R)$ on the upper half
plane is transitive,  we can assume without loss of generality
that $z_1=i$ for the proof of the existence of $M$. Then the
problem becomes to prove the existence of $M$ with fixed point $i$
such that $M<z_2>=i\l$. Using the transformation $z\longrightarrow
(z-i)/(z+i)$, we are led to consider linear fractional transforms
with fixed point $0$. Let $z'_2$ be the image of $z_2$ under this
transformation, then $z'_2$ is in the unit circle centered at the
origin. Write $z'_2=|z'_2|e^{i\th}$ with $0\le \th <2\pi $. Then
the transformation $z\longrightarrow e^{-i\th}z$ will have fixed
point $0$ and takes $z'_2$ to $r=|z'_2|<1$. Transforming this into
the upper half plane we have proved the existence of $M$. To prove
the uniqueness of $\l$, we suppose there exists another $M'\in
SL_2(\Bbb R)$ such that
$$
M'<z_1>=i,\,\,\,M'<z_2>=i\l ',
$$
with $\l '\ge 1$. Then we have
$$
M'M^{-1}<i\l >=i\l ',\,\,\,M'M^{-1}<i>=i.\tag 75
$$
By the second equality of (75) we see that
$$
M'M^{-1}=\ve \pmatrix \cos\th & \sin\th \\ -\sin\th & \cos\th \endpmatrix
$$
where $\th$ is real. Thus by the first equality of (75) we get
$$
i\l \cos\th +\sin\th =\l\l'\sin\th +i\l'\cos\th .
$$
Hence if $\cos\th\ne 0$ then by the identity of the imaginary part we see that $\l=\l'$. If $\cos\th =0$ then $\sin\th \ne 0$, so $\l\l'=1$, which also leads to $\l=\l'$ since we have $\l\ge 1$ and $\l'\ge 1$.
Now, from $M<z_1>=i$ and Lemma 10 we have
$$
M=\pmatrix \cos\th_1 & \sin\th_1 \\ -\sin\th_1 & \cos\th_1 \endpmatrix
\pmatrix y_1^{-1/2} & -x_1y_1^{-1/2} \\ 0 & y_1^{1/2} \endpmatrix ,
$$
and from $M<z_2>=i\l$ and Lemma 10 we have
$$
M=\pmatrix \l^{1/2} & 0 \\ 0 & \l^{-1/2} \endpmatrix
\pmatrix \cos\th_2 & \sin\th_2 \\ -\sin\th_2 & \cos\th_2 \endpmatrix
\pmatrix y_2^{-1/2} & -x_2y_2^{-1/2} \\ 0 & y_2^{1/2} \endpmatrix ,
$$
where $\th_1$ and $\th_2$ are real numbers satisfying $0\le
\th_1,\,\,\th_2<2\pi$. By these and Lemma 8 we complete the proof
of Lemma 10.
\enddemo

\proclaim{Theorem 7} Let $Z_1=\pmatrix \tau_1 & z_1 \\ z_1 &
\tau_1 \endpmatrix$ and $Z_2=\pmatrix \tau_2 & z_2 \\ z_2 & \tau_2
\endpmatrix$ be two arbitrary points in $\hat\Bbb H_2$.
Write
$$
\align &
\tau_1+z_1=x_1+iy_1,\,\,\tau_1-z_1=x_2+iy_2,\,\,\tau_2+z_2=u_1+iv_1,\,\,\tau_2-z_2=u_2+iv_2,\\
& A=\f{y_1^2 + v_1^2 + \left(x_1 - u_1\right)^2}{y_1v_1}\ge 2,\,\,
B=\f{y_2^2 + v_2^2 + \left(x_2 - u_2\right)^2}{y_2v_2}\ge 2.
\endalign
$$
Then the distance $\rho \left(Z_1,\,Z_2\right)$ of $Z_1$ and $Z_2$
is equal to
$$
\rho
\left(Z_1,\,Z_2\right)=\left(\log^2\f{A+\sqrt{A^2-4}}{2}+\log^2\f{B+\sqrt{B^2-4}}{2}\right)^{1/2}.
$$
\endproclaim

\demo{Proof} By Theorem 4 we know that there exists an $M\in
\hat\O_2$ such that
$$
M<Z_1>=iI,\,\,M<Z_2>=i\L ,
$$
where $\L =\pmatrix \l_1 & \l_2 \\ \l_2 & \l_1 \endpmatrix \in
{\Bbb R}^{(2,2)}$ is a real matrix with $\l_1\ge \l_2+1$ and
$\l_2\ge 0$. Thus by Lemma 7 we may obtain the induced matrices
$M_1,\,M_2 \in SL_2(\Bbb R)$ from $M$ such that
$$
\f12\left(M_1<\tau_1+z_1>+M_2<\tau_1-z_1>\right)=i,\,\,\,M_1<\tau_1+z_1>-M_2<\tau_1-z_1>=0,
$$
and
$$
\f12\left(M_1<\tau_2+z_2>+M_2<\tau_2-z_2>\right)=i\l_1,\,\,\,\f{\ve}{2}\left(M_1<\tau_2+z_2>-M_2<\tau_2-z_2>\right)=i\l_2.
$$
These imply respectively
$$
M_1<\tau_1+z_1>=i,\,\,M_2<\tau_1-z_1>=i,\tag 76
$$
and
$$
M_1<\tau_2+z_2>=i\left(\l_1+\ve\l_2\right),\,\,
M_2<\tau_2-z_2>=i\left(\l_1-\ve\l_2\right).\tag 77
$$
From Lemma 10 and the first equalities of (76) and (77) we get
$$
\l_1+\ve\l_2+\f1{\l_1+\ve\l_2}=A,
$$
and from Lemma 10 and the second equalities of (76) and (77) we
get
$$
\l_1-\ve\l_2+\f1{\l_1-\ve\l_2}=B.
$$
Thus in view of $A\ge 2$, $B\ge 2$, $\l_1+\ve\l_2\ge 1$ and
$\l_1-\ve\l_2\ge 1$, we get from the above two equalities that
$$
\l_1+\ve\l_2=\f{A+\sqrt{A^2-4}}{2},\,\,\l_1-\ve\l_2=\f{B+\sqrt{B^2-4}}{2}.
$$
Therefore by Lemma 5 and Theorem 6 we get
$$
\align
\rho\left(Z_1,\,Z_2\right)&=\rho\left(iI,\,i\L\right)=\left(\log^2\left(\l_1+\l_2\right)+\log^2\left(\l_1-\l_2\right)\right)^{1/2}\\
&=\left(\log^2\f{A+\sqrt{A^2-4}}{2}+\log^2\f{B+\sqrt{B^2-4}}{2}\right)^{1/2},
\endalign
$$
as what we need. The proof of the theorem is complete.
\enddemo

Next, we are going to give the parameter equation of geodesic
connecting two arbitrary different points $Z_1$ and $Z_2$ in
$\hat\Bbb H_2$. To do this, we need the following preliminary
result.

\proclaim{Lemma 11} Let $Z_1$ and $Z_2$ be any two different fixed
points in $\hat\Bbb H_2$, let $M=\pmatrix A & B \\ C & D
\endpmatrix \in \hat\O_2$ be such that
$$
M^{-1}<Z_1>=iI,\,\,\,M^{-1}<Z_2>=i\L ,
$$
with $\L=\pmatrix \l_1 & \l_2 \\ \l_2 & \l_1 \endpmatrix$ and
$\l_1\ge \l_2 +1\ge 1$. Suppose $Z=Z(s)$ is the parameter equation
of the geodesic connecting the points $Z(0)=iI$ and $Z(s_0)=i\L$,
where the parameter $s$ with $0\le s\le s_0$ denotes the length of
the arcs of the geodesic, then
$$
Z=W(s):=M<Z(s)>
$$
is the parameter equation of geodesic connecting the points
$Z_1=M<Z(0)>=W(0)$ and $Z_2=M<Z(s_0)>=W(s_0)$, with $0\le s\le
s_0$ being also the length of the arcs of this geodesic. So $W(s)$
is irrelevant to the choice of the matrix $M$ above.
\endproclaim

\demo{Proof} By Lemma 6, we only need to prove
$$
\ddot W(s)+i\dot W(s)V(s)^{-1}\dot W(s)=0, \tag 78
$$
under the assumption that
$$
\ddot Z(s)+i\dot Z(s)Y(s)^{-1}\dot Z(s)=0, \tag 79
$$
where $V(s)={\text Im}W(s)$ and $Y(s)={\text Im}Z(s)$.
Write
$$
Z(s)=\pmatrix \tau (s) & z(s) \\ z(s) & \tau (s) \endpmatrix .
$$
Then from Lemma 7 we know that there exists $M_1,\,M_2\in
SL_2(\Bbb R)$ such that
$$
W(s)=M<Z(s)>=\pmatrix \f{M_1<\tau
(s)+z(s)>+M_2<\tau (s)-z(s)>}{2} & \f{\ve \left(M_1<\tau
(s)+z(s)>-M_2<\tau (s)-z(s)>\right)}{2}
 \\ \f{\ve \left(M_1<\tau (s)+z(s)>-M_2<\tau (s)-z(s)>\right)}{2} & \f{M_1<\tau (s)+z(s)>+M_2<\tau (s)-z(s)>}{2} \endpmatrix
$$
Thus by multiplying $p^{-1}$ and $p$ from the left and right
respectively to both sides of (79) we see that (79) becomes $$
\left\{
\aligned &\left(\tau (s)+z(s)\right)^{..}+i{\left(\tau (s)+z(s)\right)^{.}}^2{\text Im}\left(\tau (s)+z(s)\right)^{-1}=0,   \\
&\left(\tau (s)-z(s)\right)^{..}+i{\left(\tau
(s)-z(s)\right)^{.}}^2{\text Im}\left(\tau (s)-z(s)\right)^{-1}=0.
\endaligned
 \right.\tag 80
$$
And similarly (78) can be rewritten as
$$
\left\{
\aligned &\left(M_1<\tau (s)+z(s)>\right)^{..}+i{\left(M_1<\tau (s)+z(s)>\right)^{.}}^2{\text Im}\left(M_1<\tau (s)+z(s)>\right)^{-1}=0,   \\
&\left(M_2<\tau (s)-z(s)>\right)^{..}+i{\left(M_2<\tau
(s)-z(s)>\right)^{.}}^2{\text Im}\left(M_2<\tau
(s)-z(s)>\right)^{-1}=0.
\endaligned
 \right.\tag 81
$$
Then our problem becomes to prove (81) under (80). Notice that
$\tau (s)\pm z(s)$ is in the upper half plane and $M_1$ and $M_2$
are in $SL_2(\Bbb R)$. So it is sufficient to prove, for any
$M=\pmatrix a & b \\ c & d \endpmatrix \in SL_2(\Bbb R)$ and any
$z(s)$ in the upper half plane,
$$
\left(M<z(s)>\right)^{..}+i{\left(M<z(s)>\right)^{.}}^2{\text
Im}\left(M<z(s)>\right)^{-1}=0, \tag 82
$$
under the assumption that
$$
\left(z(s)\right)^{..}+i{\left(z(s)\right)^{.}}^2{\text
Im}\left(z(s)\right)^{-1}=0. \tag 83
$$
Note that, by simple computation,
$$
\align & {\text Im}M<z(s)>={\text
Im}\f{az(s)+b}{cz(s)+d}=\f{{\text Im}z(s)}{|cz(s)+d|^2},\\
&
\left(M<z(s)>\right)^{.}=\left(\f{az(s)+b}{cz(s)+d}\right)^{.}=\f{\dot
z(s)}{(cz(s)+d)^2},\\
& \left(M<z(s)>\right)^{..}=\f{(cz(s)+d)\ddot z(s)-2c\dot
z(s)^2}{(cz(s)+d)^3}.
\endalign
$$
Substituting these into (82) we can rewrite it as
$$
\left({\text
Im}z(s)\right)(cz(s)+d)\ddot z(s)-2c\left({\text
Im}z(s)\right)\dot z(s)^2+i\dot z(s)^2\left(c\bar z(s+d)\right)=0.
$$
Using (83), this can be written further as
$$
-i(cz(s)+d)\dot z(s)^2-2c\left({\text Im}z(s)\right)\dot
z(s)^2+i\dot z(s)^2\left(c\bar z(s+d)\right)=0,
$$
which follows clearly from the trivial identity
$$
-i\left(cz(s)+d\right)-2c{\text Im}z(s)+i\left(c\bar
z(s)+d\right)=0.
$$
The proof of this lemma is complete.
\enddemo

\proclaim{Theorem 8} Let $Z_1=\pmatrix \tau_1 & z_1 \\ z_1 &
\tau_1 \endpmatrix$ and $Z_2=\pmatrix \tau_2 & z_2 \\ z_2 & \tau_2
\endpmatrix$ be any two fixed points in $\hat\Bbb H_2$. Write
$$
\align &
\tau_1+z_1=x_1+iy_1,\,\,\tau_1-z_1=x_2+iy_2,\,\,\tau_2+z_2=u_1+iv_1,\,\,\tau_2-z_2=u_2+iv_2,\\
& A=\f{y_1^2 + v_1^2 + \left(x_1 - u_1\right)^2}{y_1v_1}\ge 2,\,\,
B=\f{y_2^2 + v_2^2 + \left(x_2 - u_2\right)^2}{y_2v_2}\ge 2, \\
& \l =\f{A+\sqrt{A^2-4}}{2}\ge 1,\,\, \t\l
=\f{B+\sqrt{B^2-4}}{2}\ge 1.
\endalign
$$
Then the geodesic connecting $Z_1$ and $Z_2$ can be written as
$$
Z=W(s)=\pmatrix \tau(s) & z(s) \\ z(s) & \tau(s)
\endpmatrix,
$$
where
$$ \align \tau(s)&=\f{x_1+x_2}2+\f{y_1}{2}\f{\l
\left(u_1-x_1\right)\left(\l^{2s/s_0}-1\right)c(\l)+iv_1\l^{s/s_0}}{\left(\l
y_1-v_1\right)\left(\l^{2s/s_0}-1\right)c(\l)+v_1}\\
&\quad +\f{y_2}{2}\f{{\t\l}
\left(u_2-x_2\right)\left({\t\l}^{2s/s_0}-1\right)c({\t\l})+iv_2{\t\l}^{s/s_0}}{\left({\t\l}
y_2-v_2\right)\left({\t\l}^{2s/s_0}-1\right)c({\t\l})+v_2},
\endalign
$$
$$
\align z(s)&=\f{x_1-x_2}2+\f{y_1}{2}\f{\l
\left(u_1-x_1\right)\left(\l^{2s/s_0}-1\right)c(\l)+iv_1\l^{s/s_0}}{\left(\l
y_1-v_1\right)\left(\l^{2s/s_0}-1\right)c(\l)+v_1}\\
&\quad -\f{y_2}{2}\f{{\t\l}
\left(u_2-x_2\right)\left({\t\l}^{2s/s_0}-1\right)c({\t\l})+iv_2{\t\l}^{s/s_0}}{\left({\t\l}
y_2-v_2\right)\left({\t\l}^{2s/s_0}-1\right)c({\t\l})+v_2},
\endalign
$$
$c(x)$ is a function defined on the interval $[1,\, \infty )$ by
$c(x):=\left(x^2-1\right)^{-1}$ if $x>1$, and $c(x):=1$ if $x=1$,
and the parameter $s$ with $0\le s\le s_0$ denotes the length of
the arc initiated at $Z_1$ in the geodesic. So the geodesic is
uniquely determined by $Z_1$ and $Z_2$.
\endproclaim

\demo{Proof} By Theorem 4 we know that there exists an $M\in
\hat\O_2$ such that
$$
M<Z_1>=iI,\,\,M<Z_2>=i\L=i\pmatrix \l_1 & \l_2 \\ \l_2 & \l_1
\endpmatrix , \tag 84
$$
where $\L$ with $\l_1\ge \l_2 +1\ge 1$ is determined uniquely by
$Z_1$ and $Z_2$. For this $M$, by Lemma 7 we see that there exist
two matrices $M_1,\,M_2\in SL_2(\Bbb R)$ such that, for any
$Z=\pmatrix \tau & z \\ z & \tau \endpmatrix \in \hat\Bbb H_2$,
the action of $M$ on $Z$ can be written as
$$
M<Z>=\pmatrix \f{M_1<\tau +z>+M_2<\tau -z>}{2} & \f{\ve
\left(M_1<\tau +z>-M_2<\tau -z>\right)}{2}
 \\ \f{\ve \left(M_1<\tau +z>-M_2<\tau -z>\right)}{2} & \f{M_1<\tau +z>+M_2<\tau -z>}{2}
 \endpmatrix . \tag 85
$$
By (85) with $Z=Z_1$, together with the first equality in (84), we
have
$$
\f12\left({M_1<\tau_1 +z_1>+M_2<\tau_1 -z_1>}\right)=i,
\,\, \f{\ve}{2} \left(M_1<\tau_1 +z_1>-M_2<\tau_1 -z_1>\right)=0,
$$
whence
$$
M_1<\tau_1 +z_1>=i,\,\,M_2<\tau_1 -z_1>=i. \tag 86
$$
Similarly, by (85) with $Z=Z_2$, together with the second equality
in (84), we have
$$
\f12\left({M_1<\tau_2 +z_2>+M_2<\tau_2 -z_2>}\right)=i\l_1, \,\,
\f{\ve}{2} \left(M_1<\tau_2 +z_2>-M_2<\tau_2 -z_2>\right)=i\l_2,
$$
whence
$$
M_1<\tau_2 +z_2>=i\left(\l_1+\ve \l_2\right)=i\l,\,\,M_2<\tau_2
-z_2>=i\left(\l_1-\ve \l_2\right)=i\t\l . \tag 87
$$
From the first equalities of (86) and (87), by Lemma 10 we have $
\l +\f1{\l }=A.$ This together with $\l\ge 1$ gives
$$
\l =\f{A+\sqrt{A^2-4}}{2}. \tag 88
$$
Similarly, from the second equalities of (86) and (87), by Lemma
10 and in view of $\t\l \ge 1$ we have
$$
\t\l =\f{B+\sqrt{B^2-4}}{2}. \tag 89
$$
Again, from Lemma 9 and the first equality of (86), we can write
$$
M_1=\pmatrix \cos\th_1 & \sin\th_1 \\ -\sin\th_1 & \cos\th_1
\endpmatrix
\pmatrix y_1^{-1/2} & -x_1y_1^{-1/2} \\
0 & y_1^{1/2}
\endpmatrix ,\tag 90
$$
where $\th_1$ is a real number with $0\le \th_1 <2\pi$.
Substituting this into the first equality of (87) we get
$$
\pmatrix y_1^{-1/2} & -x_1y_1^{-1/2} \\
0 & y_1^{1/2}
\endpmatrix <\tau_2+z_2>
=\pmatrix \cos\th_1 & -\sin\th_1 \\
\sin\th_1 & \cos\th_1
\endpmatrix <i\l>,
$$
that is
$$
\f{-\sin\th_1 +i\l\cos\th_1}{\cos\th_1 +i\l\sin\th_1}
=\frac{u_1-x_1+iv_1}{y_1},
$$
Thus, by comparing the real and imaginary parts,
$$
\f{(\l^2-1)\sin\th_1
\cos\th_1}{(\l^2-1)\sin^2\th_1+1}=\frac{u_1-x_1}{y_1},\,\,\,
\f{\l}{(\l^2-1)\sin^2\th_1+1}=\frac{v_1}{y_1} .
$$
This clearly implies
$$
\left(\l^2-1\right)\sin^2\th_1=\frac{\l y_1-v_1}{v_1},\,\,\,
\left(\l^2-1\right)\sin\th_1\cos\th_1=\frac{\l(u_1-x_1)}{v_1}.
\tag 91
$$
Similarly, from Lemma 9 and the second equality of (86), we can
write
$$
M_2=\pmatrix \cos\th_2 & \sin\th_2 \\ -\sin\th_2 & \cos\th_2
\endpmatrix
\pmatrix y_2^{-1/2} & -x_2y_2^{-1/2}\\
0 & y_2^{1/2} \endpmatrix,\tag 92
$$
where $\th_2$ is again a real number with $0\le \th_2 <2\pi$.
Substituting this into the second equality of (87) we get
$$
\pmatrix y_2^{-1/2} & -x_2y_2^{-1/2}\\
0 & y_2^{1/2} \endpmatrix <\tau_2-z_2>
=\pmatrix \cos\th_2 & -\sin\th_2 \\
\sin\th_2 & \cos\th_2
\endpmatrix <i\tilde{\l}>,
$$
that is
$$
\f{-\sin\th_2 +i\tilde{\l}\cos\th_2}{\cos\th_2
+i\tilde{\l}\sin\th_2} =\frac{u_2-x_2+iv_2}{y_2},
$$
Thus, by comparing the real and imaginary parts,
$$
\f{\left(\tilde{\l}^2-1\right)\sin\th_2
\cos\th_2}{\left(\tilde{\l}^2-1\right)\sin^2\th_2+1}=\frac{u_2-x_2}{y_2},\,\,\,
\f{\tilde{\l}}{\left(\tilde{\l}^2-1\right)\sin^2\th_2+1}=\frac{v_2}{y_2}
.
$$
This clearly implies
$$
\left(\tilde{\l}^2-1\right)\sin^2\th_2=\frac{\t{\l}
y_2-v_2}{v_2},\,\,\,
\left(\t{\l}^2-1\right)\sin\th_2\cos\th_2=\frac{\t{\l}(u_2-x_2)}{v_2}.
\tag 93
$$
On the other hand, by Lemma 6 it is known that the geodesic
connecting the points $iI$ and $i\L$ can be expressed as
$$
Z=Z(s)=i\pmatrix\frac{(\lambda_1+\lambda_2)^{s/s_0}+(\lambda_1-\lambda_2)^{s/s_0}}{2}
& \frac{(\lambda_1
+\lambda_2)^{s/s_0}-(\lambda_1-\lambda_2)^{s/s_0}}{2}\\
\frac{(\lambda_1+\lambda_2)^{s/s_0}-(\lambda_1-\lambda_2)^{s/s_0}}{2}
&
\frac{(\lambda_1+\lambda_2)^{s/s_0}+(\lambda_1-\lambda_2)^{s/s_0}}{2}\endpmatrix
,
$$
where $s$ with $0\le s\le s_0$ is the length of the arc of the
geodesic. Hence by Lemma 11, the parameter equation of the
geodesic connecting the points $Z_1$ and $Z_2$ can be written as
$$
Z=W(s):=M^{-1}<Z(s)>, \tag 94
$$
where $s$ with $0\le s\le s_0$ is also the length of the arcs of
the geodesic. Write
$$
W(s)=M^{-1}<Z(s)>=\pmatrix \tau(s) & z(s) \\ z(s) & \tau(s)
\endpmatrix ,
$$
then we have $M<W(s)>=Z(s)$. Hence by (85) we obtain
$$
\align &
\f12\left(M_1<\tau(s)+z(s)>+M_2<\tau(s)-z(s)>\right)=\f{i}2\left((\lambda_1+\lambda_2)^{s/s_0}+(\lambda_1-\lambda_2)^{s/s_0}\right),\\
&
\f{\ve}2\left(M_1<\tau(s)+z(s)>-M_2<\tau(s)-z(s)>\right)=\f{i}2\left((\lambda_1+\lambda_2)^{s/s_0}-(\lambda_1-\lambda_2)^{s/s_0}\right).
\endalign
$$
Thus
$$
\aligned & \quad M_1<\tau(s)+z(s)>
\\
&
=\f{i}2\left(\left((\lambda_1+\lambda_2)^{s/s_0}+(\lambda_1-\lambda_2)^{s/s_0}\right)+\ve
\left((\lambda_1+\lambda_2)^{s/s_0}-(\lambda_1-\lambda_2)^{s/s_0}\right)\right)\\
&=i\l^{s/s_0},
\endaligned\tag 95
$$
and
$$
\aligned &\quad
 M_2<\tau(s)-z(s)> \\
 &=\f{i}2\left(\left((\lambda_1+\lambda_2)^{s/s_0}+(\lambda_1-\lambda_2)^{s/s_0}\right)-\ve
\left((\lambda_1+\lambda_2)^{s/s_0}-(\lambda_1-\lambda_2)^{s/s_0}\right)\right)\\
&=i\t{\l}^{s/s_0}.
\endaligned \tag 96
$$
Substituting (90) into (95), we obtain
$$
\align
&\pmatrix y_1^{-1/2} & -x_1y_1^{-1/2}\\
0 & y_1^{1/2}\endpmatrix <\tau(s)+z(s)> =\pmatrix \cos\th_1 &
-\sin\th_1 \\ \sin\th_1 & \cos\th_1
\endpmatrix <i\l^{s/s_0}>\\
&=\f{i\l^{s/s_0}\cos\th_1-\sin\th_1}{i\l^{s/s_0}\sin\th_1+\cos\th_1}
=\f{\left(\l^{2s/s_0}-1\right)\sin\th_1\cos\th_1+i\l^{s/s_0}}
{\left(\l^{2s/s_0}-1\right)\sin^2\th_1+1}.
\endalign
$$
Hence
$$
\tau(s)+z(s)=x_1+y_1\f{\left(\l^{2s/s_0}-1\right)\sin\th_1\cos\th_1+i\l^{s/s_0}}
{\left(\l^{2s/s_0}-1\right)\sin^2\th_1+1}.
$$
Thus, if $\l=1$, then
$$
\tau(s)+z(s)=x_1+iy_1\l^{s/s_0}.
$$
If $\l\neq 1$, then by (91),
$$
\tau(s)+z(s)=x_1+y_1\f{\left(\l^{2s/s_0}-1\right)\f{\l\left(u_1-x_1\right)}{\left(\l^2-1\right)v_1}+i\l^{s/s_0}}
{\left(\l^{2s/s_0}-1\right)\f{\l
y_1-v_1}{\left(\l^2-1\right)v_1}+1}.
$$
Gathering together these two cases, we get
$$
\tau(s)+z(s)=x_1+y_1\f{\l\left(u_1-x_1\right)c(\l)\left(\l^{2s/s_0}-1\right)+iv_1\l^{s/s_0}}
{\left(\l y_1-v_1\right)c(\l)\left(\l^{2s/s_0}-1\right)+v_1}.\tag
97
$$
Similarly, substituting (92) into (96), we obtain
$$
\tau(s)-z(s)=x_2+y_2\f{\left(\t{\l}^{2s/s_0}-1\right)\sin\th_2\cos\th_2+i\t{\l}^{s/s_0}}
{\left(\t{\l}^{2s/s_0}-1\right)\sin^2\th_2+1}.
$$
Thus, if $\t{\l}=1$, then
$$
\tau(s)-z(s)=x_2+iy_2\t{\l}^{s/s_0}.
$$
If $\t{\l}\neq 1$, then by (93),
$$
\tau(s)-z(s)=x_2+y_2\f{\left(\t{\l}^{2s/s_0}-1\right)\f{\t{\l}\left(u_2-x_2\right)}{\left(\t{\l}^2-1\right)v_2}+i\t{\l}^{s/s_0}}
{\left(\t{\l}^{2s/s_0}-1\right)\f{\t{\l}
y_2-v_2}{\left(\t{\l}^2-1\right)v_2}+1}.
$$
Gathering together these two cases, we get
$$
\tau(s)-z(s)=x_2+y_2\f{\t{\l}\left(u_2-x_2\right)c(\t{\l})\left(\t{\l}^{2s/s_0}-1\right)+iv_2\t{\l}^{s/s_0}}
{\left(\t{\l}
y_2-v_2\right)c(\t{\l})\left(\t{\l}^{2s/s_0}-1\right)+v_2}.\tag 98
$$
The combination of (97) and (98) clearly implies the desired
result. The proof of Theorem 8 is complete.

\enddemo

\head 8. The symplectic measure related to the symplectic metric
\endhead

In this last section we will give an explicit formulation of the
symplectic measure ${\text d}v$ induced from the symplectic metric
${\text d}s$ given by the previous section, in terms of the usual
Euclidean measure. Our main result is the following

\proclaim{Theorem 9} As for the symplectic metric ${\text d}s$ of
$\hat\Bbb H_2$, the corresponding symplectic measure ${\text d}v$,
i.e., the volume element at a point $Z=X+iY=\pmatrix x_1+iy_1 &
x_2+iy_2
\\ x_2+iy_2 & x_1+iy_1
\endpmatrix$ in $\hat\Bbb H_2$, can be expressed in terms of the Euclidean measure element ${\text d}x_1{\text
d}x_2{\text d}y_1{\text d}y_2$ as
$$
{\text d}v =\f{4}{(y_1+y_2)^2(y_1-y_2)^2}{\text d}x_1{\text
d}x_2{\text d}y_1{\text d}y_2,
$$
which is invariant under the actions of the elements of
$\hat\O_2$.
\endproclaim

\demo{Proof} First of all, we recall a general assertion from
Riemannian geometry: For any metric ${\text d}s^2=\left({\text
d}x_1,\dots ,{\text d}x_n \right)A\ ^t\left({\text d}x_1,\dots
,{\text d}x_n \right)$ in an $n$ dimensional space, the
corresponding volume element ${\text d}v$ can be expressed as
$$
{\text d}v=\left(\det A\right)^{1/2}{\text d}x_1\dots {\text
d}x_n, \tag 99
$$
where $A$ is an $n\times n$ positive definite matrix in $\Bbb
R^{(n,n)}$. In the present situation, we have
$$
\align {\text d}s^2&={\text tr}\left(Y^{-1} {\text d}Z Y^{-1}
{\text d}\bar{Z}\right)={\text tr}\left(Y^{-1} {\text d}X Y^{-1}
{\text d}X+Y^{-1} {\text d}Y Y^{-1} {\text d}Y\right)\\
&={\text tr}\left(Y^{-1} {\text d}X Y^{-1} {\text
d}X\right)+{\text tr}\left(Y^{-1} {\text d}Y Y^{-1} {\text
d}Y\right).
\endalign
$$
Thus, if we write
$$
{\text tr}\left(Y^{-1} {\text d}X Y^{-1} {\text
d}X\right):=\left({\text d}x_1,\ {\text d}x_2\right)B\ ^t
\left({\text d}x_1,\ {\text d}x_2\right) ,
$$
and
$${\text tr}\left(Y^{-1} {\text d}Y Y^{-1}
{\text d}Y\right):=\left({\text d}y_1,\ {\text d}y_2\right)C\
^t\left({\text d}y_1,\ {\text d}y_2\right) ,
$$
with $B,\,C\in \Bbb R^{(2,2)}$, then we have
$$
\align {\text d}s^2&=\left({\text d}x_1,\ {\text d}x_2\right)B\
^t\left({\text d}x_1,\ {\text d}x_2\right)
+\left({\text d}y_1,\ {\text d}y_2\right)C\ ^t\left({\text d}y_1,\ {\text d}y_2\right)\\
&=\left({\text d}x_1,\ {\text d}x_2,\ {\text d}y_1,\ {\text
d}y_2\right)\pmatrix B & 0\\ 0 & C\endpmatrix  \
^t\left({\text d}x_1,\ {\text d}x_2,\ {\text d}y_1,\ {\text d}y_2\right)\\
&=\left({\text d}x_1,\ {\text d}x_2,\ {\text d}y_1,\ {\text
d}y_2\right)A\ ^t\left({\text d}x_1,\ {\text d}x_2,\ {\text
d}y_1,\ {\text d}y_2\right) ,
\endalign
$$
where $A=\pmatrix B & 0\\ 0 & C\endpmatrix$. Next, simple
computation yields
$$
\align
{\text tr}\left(Y^{-1} {\text d}X Y^{-1} {\text d}X\right)&={\text tr}\left(p^{-1}Y^{-1}pp^{-1}{\text d}Xpp^{-1}Y^{-1}pp^{-1}{\text d}Xp\right)\\
&={\text tr}\left((p^{-1}Yp)^{-1}p^{-1}{\text d}Xp(p^{-1}Yp)^{-1}p^{-1}{\text d}Xp\right)\\
&=\f{\left({\text d}x_1+{\text
d}x_2\right)^2}{\left(y_1+y_2\right)^2}+\f{\left({\text
d}x_1-{\text d}x_2\right)^2}{\left(y_1-y_2\right)^2},
\endalign
$$
and similarly
$$
{\text tr}\left(Y^{-1} {\text d}Y Y^{-1} {\text
d}Y\right)=\f{\left({\text d}y_1+{\text
d}y_2\right)^2}{\left(y_1+y_2\right)^2}+\f{\left({\text
d}y_1-{\text d}y_2\right)^2}{\left(y_1-y_2\right)^2}.
$$
Hence we can obtain
$$
B=C=\pmatrix(y_1+y_2)^{-2}+(y_1-y_2)^{-2} & (y_1+y_2)^{-2}-(y_1-y_2)^{-2}\\
(y_1+y_2)^{-2}-(y_1-y_2)^{-2} &
(y_1+y_2)^{-2}+(y_1-y_2)^{-2}\endpmatrix .
$$
And thus by (99) we obtain
$$ {\text d}v=(\det A)^{1/2}{\text d}x_1{\text d}x_2{\text d}y_1{\text d}y_2
=\f{4}{\left(y_1+y_2\right)^2\left(y_1-y_2\right)^2}{\text
d}x_1{\text d}x_2{\text d}y_1{\text d}y_2 ,
$$
which is what we need. The following is devoted to giving a direct
proof of the invariance of ${\text d}v$ under the actions of the
elements of $\hat{\O}_2$. Take an $M=\pmatrix A & B\\C &
D\endpmatrix\in\hat{\O}_2$ with $MQ=\ve QM$. Note that we have
$$
A=\pmatrix a_1 & a_2\\ \ve a_2 &  \ve a_1\endpmatrix,\ B=\pmatrix
b_1 & b_2\\ \ve b_2 &  \ve b_1\endpmatrix,\ C=\pmatrix c_1 & c_2\\
\ve c_2 &  \ve c_1\endpmatrix,\ D=\pmatrix d_1 & d_2\\ \ve d_2 &
 \ve d_1\endpmatrix.
$$
Write
$$
Z=\pmatrix z_1 & z_2\\z_2 & z_1\endpmatrix,\ W=\pmatrix w_1 & w_2\\w_2 & w_1\endpmatrix
$$
with $z_{j}=x_{j}+y_{j}$ and $w_{j}=u_{j}+v_{j}$ for $j=1,2$, and
let
$$
W=M<Z>=(AZ+B)(CZ+D)^{-1}.\tag 100
$$
Then our problem becomes to prove
$$
\f{4}{(y_1+y_2)^2(y_1-y_2)^2}dx_1dx_2dy_1dy_2=\f{4}{(v_1+v_2)^2(v_1-v_2)^2}du_1du_2dv_1dv_2
,
$$
which is clearly equivalent to proving
$$
\f{\partial(u_1,u_2,v_1,v_2)}{\partial(x_1,x_2,y_1,y_2)}=\f{(v_1+v_2)^2(v_1-v_2)^2}{(y_1+y_2)^2(y_1-y_2)^2}.
\tag 101
$$
Let
$$
u_1+u_2=\t{u}_1,\ u_1-u_2=\t{u}_2,\ v_1+v_2=\t{v}_1,\ v_1-v_2=\t{v}_2
$$
and
$$
x_1+\ve x_2=\t{x}_1,\ x_1-\ve x_2=\t{x}_2,\ y_1+\ve y_2=\t{y}_1,\
y_1-\ve y_2=\t{y}_2.
$$
Then
$$
\f{\partial(\t{u}_1,\t{u}_2,\t{v}_1,\t{v}_2)}{\partial(u_1,u_2,v_1,v_2)}
=\det\pmatrix 1 & 1 & 0 & 0\\ 1 & -1 & 0 & 0\\ 0 & 0 & 1 & 1\\ 0 &
0 & 1 & -1 \endpmatrix =4 ,
$$
and
$$
\f{\partial(\t{x}_1,\t{x}_2,\t{y}_1,\t{y}_2)}{\partial(x_1,x_2,y_1,y_2)}
=\det\pmatrix 1 & \ve & 0 & 0\\ 1 & -\ve & 0 & 0\\ 0 & 0 & 1 & \ve \\
0 & 0 & 1 & -\ve \endpmatrix =4.
$$
So (101) now becomes
$$
\f{\partial(\t{u}_1,\t{u}_2,\t{v}_1,\t{v}_2)}{\partial(\t{x}_1,\t{x}_2,\t{y}_1,\t{y}_2)}=\f{\t{v}_1^2\t{v}_2^2}{\t{y}_1^2\t{y}_2^2}.
\tag 102
$$
Using $p^{-1}$ and $p$ to multiply both sides of (100) from the
left and the right respectively, we can obtain by direct
computations,
$$
\pmatrix w_1+w_2 & 0\\ 0 & w_1-w_2\endpmatrix =\pmatrix
\f{(a_1+\ve a_2)(z_1+\ve z_2)+(b_1+\ve b_2)}{(c_1+\ve c_2)(z_1+\ve
z_2)+(d_1+\ve d_2)} & 0\\ 0 & \f{(a_1-\ve a_2)(z_1-\ve
z_2)+(b_1-\ve b_2)}{(c_1-\ve c_2)(z_1-\ve z_2)+(d_1-\ve d_2)}
\endpmatrix .
$$
This implies for $j=1,2$,
$$
\t{u}_j+i\t{v}_j=\f{\t{a}_j(\t{x}_j+i\t{y}_j)+\t{b}_j}{\t{c}_j(\t{x}_j+i\t{y}_j)+\t{d}_j},\tag
103
$$
where
$$
\align
&\t{a}_1=a_1+\ve a_2,\ \t{b}_1=b_1+\ve b_2,\ \t{c}_1=c_1+\ve c_2,\ \t{d}_1=d_1+\ve d_2,\\
&\t{a}_2=a_1-\ve a_2,\ \t{b}_2=b_1-\ve b_2,\ \t{c}_2=c_1-\ve c_2,\
\t{d}_2=d_1-\ve d_2.
\endalign
$$
Thus we have
$$
\f{\partial\t{u}_j}{\partial\t{x}_j}+i\f{\partial\t{v}_j}{\partial\t{x}_j}
=A_j,\,\,
\f{\partial\t{u}_j}{\partial\t{y}_j}+i\f{\partial\t{v}_j}{\partial\t{y}_j}=iA_j
,
$$
where, for $j=1,\,2$,
$$
A_j=\f{\t{a}_j\t{d}_j-\t{b}_j\t{c}_j}{(\t{c}_j(\t{x}_j+i\t{y}_j)+\t{d}_j)^2}.
$$
Hence
$$
\aligned
\f{\partial(\t{u}_1,\t{v}_1,\t{u}_2,\t{v}_2)}{\partial(\t{x}_1,\t{y}_1,\t{x}_2,\t{y}_2)}
&=\det\pmatrix {\text Re} A_1 & -{\text Im} A_1 &0 & 0\\ {\text
Im} A_1 & {\text Re} A_1 & 0 & 0\\ 0 & 0 & {\text Re} A_2
& -{\text Im} A_2\\ 0 & 0 & {\text Im} A_2 & {\text Re} A_2 \endpmatrix\\
&=\left(({\text Re} A_1)^2+({\text Im} A_1)^2\right)\left(({\text Re} A_2)^2+({\text Im} A_2)^2\right)\\
&=|A_1|^2|A_2|^2 .
\endaligned\tag 104
$$
Note that for any linear fractional transformation
$g<z>=\f{az+b}{cz+d}$, there holds
$$
{\text Im} g<z>=\f{ad-bc}{|cz+d|^2}{\text Im} z.
$$
Using this equality with $\pmatrix a & b \\ c & d \endpmatrix =
\pmatrix \t a_j & \t b_j \\ \t c_j & \t d_j \endpmatrix$ to (103)
with $z=\t{x}_j+i\t{y}_j$, and in view of $\pmatrix \t a_j & \t
b_j \\ \t c_j & \t d_j \endpmatrix \in SL_2(\Bbb R)$ by Lemma 8,
we get
$$
|A_j|^2=\f{\t{v}_j^2}{\t{y}_j^2}\f{\left(\t{a}_j\t{d}_j-\t{b}_j\t{c}_j\right)^2}{\left(\t{a}_j\t{d}_j-\t{b}_j\t{c}_j\right)^2}=\f{\t{v}_j^2}{\t{y}_j^2}.
$$
Therefore by (104) we obtain
$$
\f{\partial(\t{u}_1,\t{v}_1,\t{u}_2,\t{v}_2)}{\partial(\t{x}_1,\t{y}_1,\t{x}_2,\t{y}_2)}
=\f{\t{v}_1^2\t{v}_2^2}{\t{y}_1^2\t{y}_2^2},
$$
as desired by (102), i.e., (101). The proof of Theorem 9 is
complete.
\enddemo


\mpb \centerline {\bf References}

\vskip 0.5cm



\re{[AZ]} A.N. Andrianov and V.G. Zhuravlev, {\it Modular forms
and Hecke operators,} Translated from the 1990 Russian original by
Neal Koblitz, Translations of Mathematical Monographs, 145(1995),
AMS, Providence, RI.

\re{[B]} S. B{\"{o}}cherer, {\it Siegel modular forms and theta
series,} Proc. Symp. Pure Math., 49(1989), Part 2, pp. 3-17.


\re{[C]} E. Cartan, {\it Sur les domaines born{\'{e}}s
homog{\`{e}}nes de l'espace de $n$ variables complexes,}
Abhandlungen aus dem Mathematischen Seminar der Hansischen
Universit{\"{a}}t, 11(1936), pp.116-162.

\re{[CP]} M. Courtieu and A. Panchishkin, {\it Non-Archimedean
L-functions and arithmetical Siegel modular forms, Second
edition,} Lecture Notes in Mathematics, 1471, Springer-Verlag,
Berlin, 2004.

\re{[DI]} W. Duke and {\"{O}}. Imamo{\v {g}}lu, {\it Siegel
modular forms of small weight,} Math. Annalen, 310 (1998), pp.
73-82.

\re{[Du]} N. Dummigan, {\it Period ratios of modular forms,} Math.
Ann., 318 (2000), pp.621-636.

\re{[EZ]} M. Eichler and D. Zagier, {\it The theory of Jacobi
forms,} Progress in Mathematics, 55(1985), Birkhuser Boston, Inc.,
Boston, MA.


\re{[FV]} C. Faber and G. van der Geer, {\it Sur la cohomologie
des syst{\`{e}}mes locaux sur les espaces de modules des courbes
de genre 2 et des surfaces ab{\'{e}}liennes,} I, II C. R. Math.
Acad. Sci. Paris, 338(2004), No.5, pp.381-384 and No.6,
pp.467-470.

\re{[H1]} L.K. Hua, {\it Harmonic analysis of functions of several
complex variables in the classical domains,} Transl. Math.
Monographs, Vol.6(1963), AMS.

\re{[H2]} L.K. Hua and I. Reiner, {\it On the generators of the
symplectic modular group,} Transact. AMS 65 (1949), pp.415-426.





\re{[IS]} T. Ibukiyama and N.-P. Skoruppa, {\it A vanishing
theorem for Siegel modular forms of weight one,} Abh. Math. Sem.
Univ. Hamburg, 77(2007), pp.229-235.





\re{[KS]} W. Kohnen and N.-P. Skoruppa, {\it A certain Dirichlet
series attached to Siegel modular forms of degree two,} Invent.
Math., 95(1989), No.3, pp.541-558.

\re{[M]} H. Maass, {\it Siegel's modular forms and Dirichlet
series,} Lectures Notes in Mathematics, 216(1971),
Springer-Verlag, Berlin.

\re{[P]} H. Poincare, {\it Memoire sur les fonctions Fuchsiennes,}
Acta Math., 1 (1883), pp.193-294.

\re{[RSF]} C. Ryan, Nathan, N.-P. Skoruppa and S. Fredrik, {\it
Numerical  computation of a certain Dirichlet series attached to
Siegel modular forms of degree two,} Math. Comp., 81(2012),
No.280, pp.2361-2376.

\re{[S1]} C.L. Siegel, {\it Einheiten quadratischer Formen,} Abh.
Math. Sem. Hans. Univ., 13 (1940), pp.209-239.

\re{[S2]} C.L. Siegel, {\it Symplectic geometry,} Amer. J. Math.,
65 (1943), pp.1-86.

\re{[S3]} C.L. Siegel, {\it Zur Theorie der Modulfunktionen n-ten
Grades,} Comm. Pure Appl. Math., 8 (1955), pp.677-681.

\re{[Sk1]} N.-P. Skoruppa, {\it Explicit formulas for the Fourier
coefficients of Jacobi and elliptic modular forms,} Invent. Math.,
102(1990), No.3, pp.501-520.

\re{[Sk2]} N.-P. Skoruppa, {\it Heegner cycles, modular forms and Jacobi forms,}  S$\acute{e}$m. Th$\acute{e}$or. Nombres Bordeaux(2), 3(1991), No.1, pp.93-116.

\re{[Sk3]} N.-P. Skoruppa, {\it Computations of Siegel modular forms of genus two,} Math. Comp., 58 (1992), No.197, pp.381-398.

\re{[Sk4]} N.-P. Skoruppa, {\it Jacobi forms of critical weight
and Weil representations,} Modular  forms on Schiermonnikoog,
pp.239-266, Cambridge Univ. Press, Cambridge, 2008.

\re{[SW]} N.-P. Skoruppa and E. Wolfgang, {\it SL(2,Z)-invariant spaces spanned by modular units,} Automorphic forms and zeta functions, pp.365-388, World Sci. Publ., Hackensack, NJ, 2006.

\re{[SZ1]} N.-P. Skoruppa and D. Zagier, {\it Jacobi forms and a certain space of modular forms,}
Invent. Math., 94(1988), No.1, pp.113-146.

\re{[SZ2]} N.-P. Skoruppa and D. Zagier, {\it A trace formula for Jacobi forms,} J. Reine Angew.
Math., 393(1989), pp.168-198.


\re{[WS]} E. Wolfgang and N.-P. Skoruppa, {\it Modular invariance and uniqueness of conformal characters,} Comm. Math. Phys., 174(1995), No.1, pp.117-136.

\mpb

\spb \noindent \settabs 2 \columns \+ Tianqin Wang & \cr \+ School
of Information Engineering & \cr \+ North China University of
Water Resources and Electric Power & \cr \+ Zhengzhou 450045,
Henan, P.R.China & \cr \+ Email: wangtq \@ amss.ac.cn  & \cr \+  &
\cr \+

Tianze Wang & \cr \+ School of Mathematics and Information
Sciences & \cr \+ North China University of Water Resources and
Electric Power & \cr \+ Zhengzhou 450045, Henan, P.R.China & \cr
\+ Email: wtz \@ ncwu.edu.cn  & \cr \+  & \cr \+

Hongwen Lu & \cr \+ School of Mathematics and Information Sciences
& \cr \+ North China University of Water Resources and Electric
Power & \cr \+ Zhengzhou 450045, Henan, P.R.China & \cr \+ Email:
lu-hongwen \@ 163.com  & \cr

\enddocument